\documentclass[a4paper,12pt,twoside]{amsart}
\usepackage{amsmath,amsthm,amsfonts}
\usepackage[T1]{fontenc}
\usepackage{graphics,graphicx}
\usepackage{graphicx}

\newtheorem{thm}{Theorem}[section]

\newtheorem{cor}[thm]{Corollary}
\newtheorem{lem}[thm]{Lemma}
\newtheorem{remark}[thm]{Remark}
\newtheorem{Lemma}[thm]{Lemma}
\newtheorem{prop}[thm]{Proposition}
\newtheorem{Prop}[thm]{Proposition}
\newtheorem{proposition}[thm]{Proposition}

\theoremstyle{definition}

\newtheorem{rem}[thm]{Remark}

\newcommand{\R}{{\mathbb{R}}}

\newcommand{\Q}{{\mathbb{Q}}}

\newcommand{\Z}{{\mathbb{Z}}}
\newcommand{\N}{{\mathbb{N}}}
\newcommand{\T}{{\mathbb{T}}}

\newcommand{\Cal}{\mathcal}
\newcommand{\cal}{\mathcal}

\def\M#1{{\mathbb #1}}

\def \mod {{\rm \ mod \,} }

\def \Var {{\rm Var}}

\def\l{{\underline l}} 
 
\def\0{{\underline 0}}  \def\1{{\underline 1}}

\def \and{\text{ and }}

\def \stm0{{\setminus \{\0\}}}

\def \eop{\qed}
\def \Proof {\vskip -3mm {{\it Proof}. }}
\def \proof {\vskip -3mm {{\it Proof}. }}
\def \Proof {\vskip 0mm {{\it Proof}. \ }}

\def \Var  {{\rm Var }}

\newcommand{\xbm}{(X,{\cal B},\mu)}
\newcommand{\zdr}{(Z,{\cal D},\rho)}
\newcommand{\ct}{{\cal T}}
\newcommand{\cb}{{\cal B}}
\newcommand{\ot}{\otimes}
\newcommand{\la}{\lambda}
\newcommand{\beq}{\begin{equation}}
\newcommand{\eeq}{\end{equation}}
\newcommand{\tend}[3][]{\xrightarrow[#2\to#3]{#1}}
\date{\today}

\title[Centralizer and liftable centralizer of special flows over rotations] {Centralizer and liftable centralizer of special flows over rotations}
\author{Jean-Pierre Conze \and Mariusz Lema\'nczyk}
\address{Jean-Pierre Conze, IRMAR, UMR CNRS 6625, Universit\'e de Rennes I
\vskip 0mm Campus de Beaulieu, 35042 Rennes Cedex, France}
\email{conze@univ-rennes1.fr}
\address{Mariusz Lema\'nczyk, Faculty of Mathematics and Computer Science, Nicolaus Copernicus University,
\vskip 0mm Chopin street 12/18, 87-100 Toru\'n, Poland}
\email{mlem@mat.umk.pl}

\thanks{Research supported by Narodowe Centrum Nauki grant  UMO-2014/15/B/ST1/03736.
Research supported by the special program in the framework of the Jean Morlet semester ``Ergodic Theory and Dynamical Systems
in their Interactions with Arithmetic and Combinatorics''.}
\begin{document}


\maketitle

\thispagestyle{empty}

\begin{abstract} The liftable centralizer for special flows over irrational rotations is studied. It is shown
that there are such flows under piecewise constant roof functions which are rigid and whose liftable centralizer is trivial.
\end{abstract}

\tableofcontents

\newpage

\newpage

\section*{\bf Introduction}

{\bf Flows and special flows}

Assume that $\zdr$ is a probability standard Borel space. In this paper, we will deal with measurable,
measure-preserving\footnote{We tacitly assume that these $\R$-actions are {\em free}, i.e.,\ for $\rho$-a.e.\ $z\in Z$, the map $t\mapsto T_tz$ is 1-1.} $\R$-actions,
i.e., with flows $\ct=(T_t)_{t\in\R}$ acting on $\zdr$ for which the map $(z,t)\mapsto T_tz$ is measurable and $\rho(T_tA)=\rho(A)$ for each $A\in\mathcal{D}$
and $t\in\R$. It follows that the unitary representation in $L^2\zdr$ corresponding to $\ct$ is strongly (equivalently, weakly) continuous, i.e.,\ the map
$t\mapsto T_tf$ is continuous for each $f\in L^2\zdr$, where $T_tf=f\circ T_t$.

We constantly assume ergodicity of flows under consideration.
According to Ambrose-Kakutani theorem \cite{Am} each such flow possesses a special representation, i.e., it can be represented as a special flow $T^f=(T^f_t)_{t\in\R}$,
where $T$ is an ergodic automorphism (often called a base) of a probability standard Borel space $\xbm$,
and $f:X\to \R^+$ is in $L^1\xbm$ ($f$ is often called a roof function).

Recall that $T^f$ acts on $(X^f,\cb^f,\mu^f)$, where $X^f=\{(x,s)\in X\times\R:\: 0\leq s<f(x)\}$ on which we consider the restriction
of product $\sigma$-algebra and product measure (which is normalized: $\mu^f(A)=(\mu\ot\la_{\R})(A)/\int_Xf\,d\mu$ for each $A\in\cb^f$).
Then, for all $t\in\R$ and $(x,s)\in X^f$, we have
$$T^f_t(x,s)=(T^nx,t+s-f^{(n)}(x)),$$
where $n\in\Z$ is unique such that   $f^{(n)}(x)\leq t+s<f^{(n+1)}(x)$. Here,
\begin{equation} \label{sumy}
f^{(n)}(x)=f(x)+f(Tx)+\ldots+f(T^{n-1}x)\text{ when }n>0,
\end{equation}
$f^{(0)}(x)=0$ and the cocycle identity $f^{(m+n)}(x)=f^{(n)}(x)+f^{(m)}(T^nx)$, true for all integers $m,n$ determines the values of $f^{(m)}$ for negative integers.

Of course, in general, a flow has many special representations
(with non-isomorphic bases). Originated by von Neumann \cite{vonNe}, it is a rather common and fruitful approach to study flows by choosing a suitable
special representation. From that point of view a lot of attention has been devoted to study special flows over irrational rotations, or,
more generally, over interval exchange transformations, as often they are natural special representations of interesting smooth, or smooth singular,
flows on surfaces, see e.g.\ \cite{Fa}, \cite{Fa-Ka}, \cite{Fr-Le1}, \cite{Fr-Le-Le}, \cite{Kh-Si}, \cite{Ko}, \cite{Ko1}, \cite{Ko2}, \cite{Ko3},
\cite{Ku}, \cite{Le00}, \cite{Sc}, \cite{Ul1}, \cite{Ul2}.

{\bf Centralizer}

A particular object of study in this paper is the {\em centralizer} of flows.
We recall that given a flow $\ct=(T_t)_{t\in\R}$ on $\zdr$, its centralizer $C(\ct)$ consists of all automorphisms $W$ of $\zdr$ commuting
with all $T_t$, $t\in\R$. When $C(\ct)=\{T_t:\:t\in\R\}$, then one says that $\ct$ has a {\em trivial centralizer}.
In general, $\{T_t:\:t\in\R\}\subset C(\ct)$ is a normal subgroup of $C(\ct)$ and the quotient group $C(\ct)/\{T_t:\:t\in\R\}$ is called the {\em essential centralizer} of $\ct$.
The essential centralizer can be quite big. Indeed, for example, it is uncountable when $\ct$ is rigid. Recall that rigidity means that for some sequence $r_n\to\infty$,
we have $T_{r_n}\to Id$ strongly in $L^2\zdr$.~\footnote{The fact that the essential centralizer is uncountable for rigid flows is folklore, see a proof of this fact,
e.g.\ in \cite{Ka-Le}, see also the proof of Proposition~\ref{pot3} below.} Prominent examples of rigid flows are given by the class
of  area-preserving smooth flows $\ct=(T_t)_{t\in\R}$ without fixed points on $\T^2$,  see \cite{Co-Fo-Si}, Chapter~16. Hence, such flows have uncountable essential centralizers.

{\bf Centralizer for special flows. Liftable centralizer}

First let us consider the continuous case: $X$ is a compact metric space, $f:X\to\R^+$ is continuous and $W\in C(T^f)$ acting on $X^f$  is also continuous,
i.e., it belongs to $C^{\rm top}(T^f)$ (note that $X^f$ has a natural metric making it a compact metric space, see Appendix).~\footnote{Clearly, under such assumptions,
the special flow $T^f$ is also a continuous flow.  When $T$ is uniquely ergodic with $\mu$ the unique $T$-invariant measure, also $T^f$ is uniquely ergodic;
so each continuous $W:X^f\to X^f$, $T^f_f\circ W=W\circ T^f_t$ for all $t\in\R$, preserves the measure $\mu^f$.}

A result from \cite{Markley} states that if $X$
is a torus and $T$ a minimal rotation on $X$, each element $W$ of $C^{\rm top}(T^f)$ comes from an $S\in C(T)$, i.e.,\ $Sx=x+\beta$ for some $\beta\in X$
and a continuous $g:X\to\R$ satisfying
\beq\label{funeq} f(Sx)-f(x)=g(Tx)-g(x)\text{ for all }x\in X.\eeq

To understand the meaning of the equation~\eqref{funeq}, consider it for the general setup of a special flow $T^f$: $T$ is an ergodic automorphism
of $\xbm$, $S$ is in $C(T)$, $g:X\to\R$ is measurable and
\beq\label{funeqm}
f(Sx)-f(x)=g(Tx)-g(x), \text{ for }\mu-\text{a.e.}~x\in X.\eeq Note that, up to natural identification, $X^f$ is the space of orbits
$\{(T_f)^n(x,r):\: n\in\Z\}$, $(x,r)\in X\times\R$, where $T_f:X\times\R\to X\times\R$,
\begin{equation}\label{skew1}
T_f(x,r)=(Tx,r+f(x))\text{ for each }(x,r)\in X\times\R.
\end{equation}
Now, the equation~\eqref{funeqm} means that $T_f\circ S_g=S_g\circ T_f$,  where
\begin{equation}\label{skewSg}
S_g(x,r) = (Sx, r + g(x)).
\end{equation}
So $S_g$ also acts on $X^f$ which is identified with $X\times\R/\sim$ and it commutes with the quotient vertical action of $\R$ which represents the special flow
(in these ``new coordinates''), see Section~\ref{lift} for details. It follows that each measurable solution $g$ of~\eqref{funeqm}
yields an element of the centralizer of $T^f$. This part of the centralizer (which is clearly a subgroup), we will call the {\em liftable}
centralizer of the special flow $T^f$ and denote it by $C^{\rm lift}(T^f)$ (of course $\{T^f_t:\:t\in\R\}\subset C^{\rm lift}(T^f)$).

One can ask now whether the liftable centralizer is the whole centralizer of the flow under consideration. But the answer to such a question
is clearly negative. For example if the base automorphism $T$ has trivial centralizer, so must be the liftable centralizer of $T^f$
for any roof function $f$.
\footnote{If $\ct$ is an arbitrary zero entropy and  loosely Bernoulli flow then it will have a special representation $T^f$ in which
$C(T)=\{T^n:\: n\in\Z\}$ \cite{Or-Ru-We}.} Moreover, if $C(T)$ is Abelian, then $C^{\rm lift}(T^f)$ is a nilpotent group of order at most~2
(see Section~\ref{lift} for basic properties of $C^{\rm lift}(T^f)$). In fact, the essential liftable centralizer $C^{\rm lift}(T_f)/\{T^f_t:\:t\in\R\}$
is Abelian. Hence if $\ct=(T_t)_{t\in\R}$ is an ergodic loosely Bernoulli flow, see \cite{Or-Ru-We}, (on $\zdr$) whose centralizer is not
a nilpotent group of order at most~2, then we cannot represent it over an ergodic $T$ so that $C(T)$ is Abelian
(which is the case for irrational rotations) to have $C^{\rm lift}(T^f)=C(T^f)$. This kind of general nonsense type arguments shows that, even
for special flows over irrational rotations, we cannot expect that the liftable centralizer is equal to the whole centralizer when $f$ is arbitrary.
\footnote{The situation does not change if additionally $X$ is a compact metric space and we require $f$ to be continuous. Indeed, each positive $L^1$-function is cohomologous
to a positive continuous function \cite{Ko}. We emphasize that even if $f$ is continuous we look for measurable solutions of \eqref{funeqm}.}

\vskip 2mm
{\bf Liftable centralizer for special flows over  irrational rotations}

We now assume that $X=\T$, where $\T$ stands for the additive circle represented as $[0,1[$ and $Tx=R_\alpha x=x+\alpha$ (mod~1), where $\alpha\in\R$ is irrational.
Supported by the aforementioned topological result of \cite{Markley}, we may still ask whether
$C^{\rm lift}(T^f)=C(T^f)$ when  $f$ is a ``natural'' function, meaning, more adapted to the topological or differentiable structure of the circle. As proved in \cite{Fr-Le2}, it is indeed the case
whenever $f$ is piecewise smooth with non-zero sum of jumps and $Tx=x+\alpha$ with $\alpha$ of bounded partial quotients.\footnote{In fact, the essential centralizer is finite in this case \cite{Fr-Le2}.}

One can ask whether $C^{\rm lift}(T^f)=C(T^f)$ when $f$ is smooth which, by \cite{Co-Fo-Si}, Chapter~16 is the case of smooth area-preserving flows without fixed points on $\T^2$.
In this case $T^f$ is rigid and hence the essential centralizer
$C(T^f)/\{T^f_t:\:t\in\R\}$ is uncountable.  In fact, even $C^{\rm lift}(T^f)/\{T^f_t:\:t\in\R\}$ is uncountable, i.e., there is always
an uncountable set of $\beta\in\T$ for which we can indeed solve~\eqref{funeqm} (with $Sx=x+\beta$) already when $f$ is absolutely continuous
-- this result is also rather folklore, so we postpone the proof of this fact to Appendix. However, the answer to the question whether $C^{\rm lift}(T^f)=C(T^f)$ for $f$ smooth is unknown. This phenomenon: $T^f$ is rigid, the number of $\beta$ for which~\eqref{funeqm}
can be solved with $Sx=x+\beta$ is uncountable, but the answer to the above question is unknown, still persists if we consider
$f=\sum_{n=-\infty}^\infty c_n e^{2\pi inx}$ and $c_n={\rm o}(1/|n|)$, see Proposition~\ref{omale}.

{\bf Main result}

In this paper we will study a relationship between $C^{\rm lift}(T^f)$ and $C(T^f)$ in the class of step functions
(for which the Fourier coefficients are clearly of order ${\rm O}(1/|n|)$). The main result is to show that they may give rigid flows
whose liftable centralizer is trivial. More precisely, we will consider $f=f_{a,b}:\T\to\R$  (with $a,b>0$) given by
\beq\label{ab}f(x)=\left\{\begin{array}{lll}
a& \text{if}& x\in[0,1/2[,\\
b& \text{if}& x\in [1/2,1[.\end{array}\right.\eeq
Under the mild assumption $a/b\notin \Q+\Q\alpha$ (which we assume to hold from now on), the special flow $T^f$ is weakly mixing \cite{Fr-Le-Le}, \cite{GuPa06}.

Let $\alpha\in[0,1[$ be irrational with the partial quotients $(a_n)_{n\geq1}$: $\alpha=[0; a_1, a_2,\ldots]$
and denominators $q_n$: $q_0=1$, $q_1=a_1$ and $q_{n+1}=a_{n+1}q_n+q_{n-1}$ for $n\geq1$.

Let us now state the main theorem (where statement 1a)  is taken from \cite{Ka-Le}) and \cite{Fr-Le-Le}).

\begin{thm} \label{Mainthm} Let $f=f_{a,b}$ with $a - b = 2$. For an irrational $\alpha$, let us consider the special flow $T^f$ obtained from $T = R_\alpha$
and $f=f_{a,b}$.

1) Suppose $\alpha$ has bounded partial quotients.
\hfill \break
1a) Then Ratner's property is satisfied for $T^f$ and therefore $C(T^f) $ is at most countable modulo $\{T^f_t:\:t\in\R\}$ (i.e.,\ the essential centralizer is at most countable).
In particular, $T^f$ is not rigid.
\hfill \break
1b) $C^{\rm lift}(T^f)$ is trivial modulo $\{T^f_t:\:t\in\R\}$.

2a) Suppose $\alpha$ has unbounded partial quotients. Then the special flow $T^f$ is rigid.
\hfill \break
2b) If $(q_{n_k})$ is even along a subsequence $(n_k)$ such that $a_{n_k+1} \uparrow \infty$, then this subsequence is a rigidity sequence  for $T^f$
and $C^{\rm lift}(T^f)$ is uncountable modulo $\{T^f_t:\:t\in\R\}$.
\hfill \break
2c) If there is $n_0$ such that the denominators $q_n$ of $\alpha$ are odd for $n\geq n_0$, then the functional equation
\begin{eqnarray}
f(x+\beta) - f(x) = g(x+\alpha) - g(x), \text{ for }\mu-\text{a.e }x\in\T,  \label{funeq0}
\end{eqnarray}
has no measurable solution $g:\T\to\R$ for $\beta \not \in \Z \alpha + \Z$. Equivalently, the liftable centralizer of $T^f$ is trivial: $C^{\rm lift}(T^f)=\{T^f_t:\:t\in\R\}$.
\hfill \break
2d) More generally, if there is $n_0$ such that $(a_{n_k+1})$ is bounded along the sequence of all $n_k$
such that $q_{n_k}$ is even,  then the conclusion is the same as in 2c).
\end{thm}

It follows that the flows from Theorem~\ref{Mainthm} display  a drastic change of ergodic properties
of special flows under the same roof function when changing  an irrational rotation as its base.

On one hand side, they seem to be interesting  from the point of view of recent achievements in studying Ratner's property \cite{Ra2}, \cite{Th}
in the class of  special flows over irrational rotations and interval exchange transformations:  \cite{Fa-Ka}, \cite{Fr-Le1}, \cite{Fr-Le-Le},
\cite{Ka1}, \cite{Ka2}, \cite{Ka-Ku}, \cite{Ka-Ku-Ul}.
Indeed (cf.\ Theorem~\ref{Mainthm}, case 1)), when $\alpha$ has bounded partial quotients then, as shown in \cite{Fr-Le-Le}, $T^f$ enjoys (finite) Ratner's property.
As proved recently in \cite{Ka-Le}, flows with (finite) Ratner's property have at most countable (discrete) essential centralizer, in particular such flows cannot be rigid.

On the other hand (cf.\ Theorem~\ref{Mainthm}, case 2a)), when $\alpha$ has unbounded partial quotients, $T^f$ is  rigid, hence,
it cannot possess the (finite) Ratner's property. Moreover, there are two different phenomena  which imply rigidity of $T^f$ in cases 2b) and 2c) in Theorem~\ref{Mainthm}.
To show that these phenomena are mutually exclusive, we will discuss them in Lemma~\ref{exclusive}, see also the second proof of Corollary~\ref{pot4}.

The paper is organized as follows: in Section~\ref{lift} we present some elementary properties of the liftable centralizer in a general setup.
Section~\ref{Prelim} is devoted to some reminders on cocycles and irrational rotations. In Section~\ref{sec1} we prove part~1 of Theorem~\ref{Mainthm}.
In Section~\ref{sec1a} we prove part 2a and part 2b) of Theorem~\ref{Mainthm}. In Section~\ref{sec2} we study the regularity of cocycles related to~(\ref{funeq0})
and in Section~\ref{sec3} we prove the remaining part of Theorem~\ref{Mainthm}. In Section~\ref{Except} we show the non-regularity of the relevant cocycle
for an exceptional set of values of $\beta$. Finally, in Appendix we study the centralizer for uniformly rigid flows and for smooth special flows over irrational rotations.
We also show that the essential liftable centralizer is uncountable whenever the Fourier transform of $f$ is of order ${\rm o}(1/|n|)$ and provide examples of H\"older
continuous roof functions which yield special flows with trivial liftable centralizer.

We would like to thank A. Danilenko, K.\ Fr\k{a}czek and A.\ Kanigowski for fruitful discussions on the subject.

\section{\bf Liftable centralizer of a special flow}\label{lift}
Let $T$ be an ergodic automorphism of a probability standard Borel space $\xbm$.
It is not hard to see that, up to natural identification, $(X^f,\mu^f)$ is the space of orbits $\{(T_f)^n(x,r):\: n\in\Z\}$, $(x,r)\in X\times\R$ (considered with the quotient of the product measure $\mu\ot\lambda_{\R}$),
where $T_f:X\times\R\to X\times\R$ is given by~\eqref{skew1}.
In these ``coordinates'' the special flow is the vertical action
$\sigma_t(x,r)=(x,r+t)$ on the quotient space $X\times \R/\sim$, where $\sim$ is the equivalence relation given by the partition into orbits of $T_f$.

Assume that $S\in C(T)$ and the equation~\eqref{funeqm} is satisfied for some measurable $g:X\to\R$. For the map $S_g$ defined by \eqref{skewSg},
it follows, that $T_f\circ S_g=S_g\circ T_f$, so $S_g$ also acts on $X^f$ (identified with $X\times\R/\sim$).
Moreover, $\sigma_t\circ S_g=S_g\circ\sigma_t$ for all $t\in\R$. Finally, $S_g$ determines an element  $\widetilde{S_g}\in C(T^f)$. Let
\beq\label{def:lift}
C^{\rm lift}(T^f):=\{\widetilde{S_g}:\;(S,g)\text{ satisfies }\eqref{funeqm}\}\eeq
be the {\em liftable centralizer} of $T^f$.
Note that $C^{\rm lift}(T^f)$ is a group as
\beq\label{e:lift0}
\widetilde{S_g} \circ \widetilde{R_h} = \widetilde{S_g\circ R_h}=\widetilde{(S\circ R)_{g\circ R+h}},\;
\widetilde{S_g}^{-1}=\widetilde{(S_g)^{-1}}=\widetilde{S^{-1}_{-g\circ S^{-1}}}.\eeq
Furthermore, for each $t\in\R$, $\widetilde{Id_t}=T^f_t$ as $Id_t(x,r)=\sigma_t(x,r)$ (where we identify $t$ with the constant function $x \to t$).
It follows that
\beq\label{e:lift1}
\{T^f_t:\;t\in\R\}\subset C^{\rm lift}(T^f).\eeq
By the same token, if $g$ is a solution of~\eqref{funeqm} then so is $g+t$ (and these exhaust all measurable solutions because of ergodicity of $T$).
Hence $\widetilde{S_{g+t}}=\widetilde{S_g}\circ T^f_t$, $t\in \R$. On the other hand, using \eqref{funeqm}, we have, for each $k\in\Z$,
\beq\label{e:lift2}
\widetilde{S_g}=\widetilde{S_g\circ(T_f)^k}.\eeq

\begin{Prop}\label{l:lift1} a) The equality $\widetilde{S_g}=\widetilde{R_h}$ holds if and only if $S_g=R_h \circ (T_f)^k$ for some $k\in\Z$.

b) If $\widetilde{S_g} \in C^{\rm lift}(T^f)$ satisfies $\widetilde{S_g}^s=Id$, then there exists $k\in\Z$ such that
$S^s=T^k$. In other words, a finite order liftable element of $C(T^f)$ must be  a lift of a root of a power of $T$.
Moreover, if $C(T)$ is trivial, so is $C^{\rm lift}(T^f)$.
\end{Prop}
\proof a) For $\mu\ot\lambda_{\R}$-a.e.\ $(x,r)\in X\times\R$, we have
$$S_g(\{(T_f)^n(x,r):\:n\in\Z\})=R_h(\{(T_f)^n(x,r):\:n\in\Z\}),$$
hence, with $k=k(x,r)$,
\beq\label{e:lift4}S_g(x,r)=R_h\circ (T_f)^k(x,r).\eeq
Since the number of $k$ is countable and, for a given $k$, the set of $(x,r)$ for which~\eqref{e:lift4} holds is measurable and $T_f$-invariant,
by the ergodicity of $T$, we obtain that $Sx=R\circ T^kx$ for $\mu$-a.e.\ $x\in X$, for some fixed $k\in\Z$.

b) The relation $\widetilde{S_g}^s = Id$ is equivalent to: for a.e. $(x, r)$, there is $k = k(x,r)$ such that:
$$(S^s x, r + \sum_{i=0}^{s-1} g(S^i x)) = (T^k x, r + \sum_{i=0}^{k-1} f(T^ix)).$$
As above, we obtain that this relation holds for some fixed $k$ and so $S^r = T^k$.
\eop

\begin{Prop}\label{p:lift1}
Assume that $C(T)$ is an Abelian group. Then $C^{\rm lift}(T^f)$ is a nilpotent group of order at most~2.
\end{Prop}
\proof \ Let $S,R$ be in $C(T)$ such that $S\circ R=R\circ S$ and~\eqref{funeqm} is satisfied for $(S, g)$ and $(R, h)$, respectively. Then, using~\eqref{e:lift0},
we obtain for the commutator:
\begin{eqnarray*}
&&(S_g \circ {R_h} \circ S_g^{-1} \circ R_h^{-1})(x, r) = \\
&&\ \  (x, r + g (S^{-1} x) - g(S^{-1} R^{-1} x) - h(R^{-1} x) + h (R^{-1} S^{-1} x)).
\end{eqnarray*}
Using \eqref{funeqm} for $(S, g)$ and $(R, h)$, a simple calculation shows that $g (S^{-1} x) - g(S^{-1} R^{-1} x) - h(R^{-1} x)
+ h (R^{-1} S^{-1} x)$ is $T$-invariant, hence a.e. equal to a constant $t$. It follows:
$$\widetilde{S_g}\circ \widetilde{R_h}\circ \widetilde{S_g}^{-1} \circ \widetilde{R_h}^{-1} = T^f_t$$
for some $t\in R$. Since $\{T^f_t:\: t\in\R\}$ is a  subgroup of the center of $C(T^f)$, we have proved the following result.
\eop

\begin{remark}\label{uw1} \em We would like to argue that, in general, $C^{\rm lift}(T^f)$ is neither a closed subgroup nor dense in $C(T^f)$.
For this aim, consider any ergodic, rigid and loosely Bernoulli flow $(R_t)$ on $\zdr$ (with the $\R$-action given by $(R_t)$ free) for which
\beq\label{e:lift5a}   \overline{\{R_t:\:t\in\R\}}\neq C((R_t)_{t\in\R}).\eeq
Note that rigidity is equivalent to:
\beq\label{e:lift5}   \{R_t:\: t\in\R\}\neq \overline{\{R_t:\:t\in\R\}}\eeq
(cf.\ the proof of Proposition~\ref{pot3}).
Clearly, properties~\eqref{e:lift5a} and~\eqref{e:lift5} are invariants of isomorphism. Now, take a special representation $T^f$ of the flow $(R_t)$ in which
$C(T)=\{T^n:\:n\in\Z\}$. Then $C^{\rm lift}(T^f)=\{T^f_t:\:t\in\R\}$. But   by~\eqref{e:lift5} and~\eqref{e:lift5a}, $C^{\rm lift}(T^f)$ is neither closed nor dense in $C(T^f)$.
\end{remark}

\section{\bf Preliminaries}\label{Prelim}

Let $\beta$ be a real number in $]0, 1[$. With $F := 1_{[0, \frac12[} - 1_{[\frac12, 0[}$, we consider the cocycle generated  over the rotation
$R_\alpha : x \to x + \alpha \mod 1$ by
\begin{eqnarray}
\Phi_\beta := \frac12 F(.-\beta) - \frac12 F. \label{cocyphi}
\end{eqnarray}
Equation (\ref{funeq0}) (where $f= 1_{[0, \frac12[} - 1_{[\frac12, 0[}$) reads $\Phi_\beta = R_\alpha g - g$.  One of our goals of this and the following sections is to show that under the assumptions of Theorem~\ref{Mainthm}, case 2c, on $\alpha$, $\Phi_\beta$ is not a coboundary, i.e.,\ equation $\Phi_\beta = R_\alpha g - g$
has no measurable solution $g$ if $\beta \not \in \Z\alpha + \Z$.
As a matter of fact,
we examine for $\Phi_\beta$ the following three properties of increasing strength:
\hfill \break (I) $\Phi_\beta$ is not a coboundary,
\hfill \break (II) $\Cal E (\Phi_\beta) \not = \{0\}$,
\hfill \break (III) $R_{\alpha, \Phi_\beta}$ is ergodic (as a skew product $(x,r) \to (R_\alpha x, r +  \Phi_\beta(x))$ on $\T\times\Z$).

Clearly if $\beta \in \Z\alpha + \Z$, then $\Phi_\beta$ is a coboundary. We exclude such values of $\beta$ which will be called {\it trivial}.

Observe that, if the group $\Cal E (\Phi_\beta)$ of finite essential values of $\Phi_\beta$ is not reduced to $\{0\}$,
then $\Phi_\beta$ is not a coboundary. We are going to show that, outside an exceptional set of values of $\beta$,
$\Cal E (\Phi_\beta) \not = \{0\}$, which implies that $\Phi_\beta$ is regular and is not a coboundary.
It remains an exceptional set of non trivial values for which $\Phi_\beta$ is not a coboundary, but can be non regular, hence non ergodic
(Theorem \ref{val-ess}).

To summarize, we will show:
\hfill \break $\bullet$ for every non trivial value of $\beta$, $\Phi_\beta$ is not a coboundary,
\hfill \break $\bullet$ for most of the values of $\beta$, it is regular,
\hfill \break$\bullet$  for an exceptional set of non trivial values of $\beta$, it is non regular, hence non ergodic.

Let us first recall some facts about essential values and useful tools in the study of cocycles (cf. \cite{Sc77}, see also \cite{Aa97}, \cite{CoRa09}).

\vskip 3mm
{\bf Reminders on cocycles}

Let $(\Phi^{(n)})$ be the cocycle (cf.~\eqref{sumy}) over an ergodic dynamical system $(X, \mu, T)$ generated by a measurable $\Phi:X\to G$,~\footnote{In what follows,
often, we call $\Phi$ itself a cocycle.} where $G = \Z^d$ or $\R^d$.  Denote by $T_\Phi$ the corresponding skew product map
$T_\Phi(x,g) = (T x, g + \Phi(x))$, $(x,g)\in X\times G$. An element $a \in G \cup \{\infty\}$
is called an {\it essential value} of the cocycle $(\Phi^{(n)})$ if, for every neighborhood $V(a)$ of $a$, for every measurable subset $B$ of positive measure,
\begin{eqnarray}
\mu(B\cap T^{-n} B \cap \{x\in X: \Phi^{(n)}(x) \in V(a) \}\bigr) > 0, {\rm \ for \ some \ } n \in \Z. \label{visitVa}
\end{eqnarray}

We denote by ${\overline {\cal E}}(\Phi)$ the  {\it set of essential values} of the cocycle $(\Phi^{(n)})$ and by ${\cal E}(\Phi) = {\overline {\cal E}}(\Phi)\cap
G$ the set of {\it finite} essential values.

A cocycle $\Phi$ is called a {\em coboundary}, if there exists a measurable $g:X\to G$ such that $\Phi(x)=g(Tx)-g(x)$ for $\mu$-a.e.\ $x\in X$. Two cocycles with values
in $G$ are said to be {\it cohomologous} if their difference is a coboundary.
Two cohomologous cocycles have the same set of essential values. If $\Phi$ is not a coboundary, then $\infty$ is an essential value of the cocycle generated by $\Phi$.
Hence, $\Phi$ is a coboundary if and only if ${\overline {\cal E}}(\Phi) =\{ 0 \}$.

The set ${\cal E}(\Phi)$ is a closed subgroup of $G$ which coincides with the group of periods $p$ of the measurable
$T_\Phi$-invariant functions on $X \times G$, i.e., the elements $p \in G$ such that, for every $T_\Phi$-invariant measurable $H$, we have
$H(x,y +p) = H(x,y), \, \mu \otimes m- a.e.$ ($m=m_G$ stands for a Haar measure on $G$). In particular, ${\cal E}(\Phi) = G$ if and only if $(X \times G, \mu \otimes m, T_\Phi)$ is ergodic.

The cocycle defined by $\Phi$ is {\it regular}, if $\Phi$ is cohomologous to a cocycle with values in a closed subgroup $G_0$ of $G$ and ergodic for the action on $X \times G_0$.
More explicitly, $\Phi$ is regular if there exists a measurable function $\eta : X \rightarrow G$ such that
$\Phi := \Psi + \eta - \eta \circ T$ $\mu$-a.e., $\Psi$ has its values in $G_0$ and $T_\Psi: (x,h) \rightarrow (T x, h + \Psi(x))$
is ergodic for the product measure $\mu \otimes m_{G_0}$ on $X \times G_0$. The group $G_0$ in this definition is necessarily ${\cal E}(\Phi)$.

In the regular case there is a ``nice'' ergodic decomposition of the measure $\mu \otimes m$ for the skew product map: any $T_\Phi$-invariant
function can be written as $V(y - \eta(x))$ for a function $V$ which is invariant by the translations by elements of $G_0$.
If the cocycle is non~regular, then the ergodic decomposition of $\mu \otimes m$ is based on a family of measures $\mu_x$ ($x\in X$) defined on $X$. Moreover, the measures $\mu_x$ are
infinite, singular with respect to the measure $\mu$ and there are uncountably many of them pairwise mutually singular.

A way to prove the existence of essential values is to use the following lemma:
\begin{lem} \label{supp}(\cite{LePaVo96})   If $(r_{n})$ is a rigidity sequence for $T$ and $(\Phi^{(r_n)})_{\ast} \mu \to \nu$ weakly on
$G \cup \{\infty\}$, then ${\rm supp}(\nu) \subset  \overline {\Cal E} (\Phi)$.
\end{lem}
A form of this criterium adapted to cocycles with values in $\Z$ is the following (\cite{Co09}):

If $a \in  G \cup \{\infty\}$ is such that there exist $\delta > 0$ and a rigidity sequence $(r_n)_{n \ge 1}$ for $T$ such that
$\mu(\{x\in X: \Phi^{(r_n)}(x) \in V(a) \}) \ge \delta$, for every neighborhood $V(a)$ of $a$, for $n$ large enough, then $a \in \overline {\Cal E}(\Phi)$.
If such an element $a$ exists and $\not \in \{0, \infty\}$, then $\Phi$ is not a coboundary.

In particular, if there exist $\delta > 0$ and a rigidity sequence $(r_n)_{n \ge 1}$ for $T$ such that
$\mu(\{x\in X: |\Phi^{(r_n)}(x)| \geq M \}) \ge \delta$, for every $M \geq 1$, for $n$ big enough, then $\infty$ is an essential value and
$\Phi$ is not a coboundary.

We will also use implicitly the following remarks: Let $f$ be a measurable $\Z$-valued function. Then if  $f$ is a $T$-coboundary in $\R$, it is a coboundary in $\Z$.
Moreover, if $T_f$ is ergodic for its action on $X \times \Z$, then the $T_f$-invariant functions on $X \times \R$ are the 1-periodic functions depending
only on the second coordinate.

\vskip 3mm
{\bf Reminders on continued fractions}

For $u \in \R$, $\|u\|$ denotes its distance to the integers: $\|u\|:= \inf_{n \in \Z} |u -n| = \min( \{u\}, 1 - \{u\}) \in [0, \frac12]$.
We will need the following inequalities:
\begin{eqnarray}
&&2 |x| \leq |\sin \pi x| \leq \pi |x|, \text{ for } |x| \leq \frac12, \label{sinA}\\
&&2 \|x\| \leq |\sin \pi x| \leq \pi \|x\|, \ \forall x \in \R. \label{sinB}
\end{eqnarray}

Let $\alpha\in[0,1[$ be an irrational number. Then, for each $n\geq1$, we write
$\displaystyle \alpha ={p_n \over q_n} + {\theta_n \over q_n}$, where $p_n $ and $q_n$ are the numerators and denominators of $\alpha$. Recall that
\begin{equation}
{1\over q_{n+1}+q_n} \le \|q_n \alpha\|  = |\theta_n| \le {1\over q_{n+1}}
= {1 \over a_{n+1} q_n+q_{n-1}},\label{f_3}\end{equation}
\begin{equation}
{1\over a_{n+1}+2} \leq q_n\|q_n \alpha\| = q_n |\theta_n| < {q_n\over q_{n+1}} < {1 \over a_{n+1}}, \label{f_3b}\end{equation}
\begin{equation}
\|k \alpha \| \geq \|q_{n-1} \alpha\| \geq \ {1\over q_{n}+q_{n-1}} \geq \ {1\over 2 q_{n}}, \text{ for}\  1\le k < q_{n}. \label{f_4}
\end{equation}
We have also
\begin{equation}\label{aa1}
(-1)^{n-1} p_n q_{n-1} = 1 + (-1)^{n-1} p_{n-1} q_n \end{equation}
and
\begin{eqnarray}
&&\|q_n \alpha\| = (-1)^n (q_n \alpha - p_n), \ \theta_n = (-1)^{n} \|q_n \alpha\|, \ \alpha = {p_n \over q_n}
+ (-1)^{n} {\|q_n \alpha\| \over q_n}. \label{f_5}
\end{eqnarray}

\begin{remark} \label{oddeven1}
\em If the denominators $q_n$ of $\alpha$ are odd for $n\geq n_0$, then the partial quotients $a_n$ are even for $n\geq n_0+2$.
Conversely, if the partial quotients $a_n$ of $\alpha$ are even for $n\geq n_0$, then for $n \geq n_0-1$ either all denominators are odd or are alternatively odd and even.
\end{remark}

In the proof of Theorem \ref{Mainthm} below, we will use the following lemma (\cite{KrLi91}, \cite{Co09}):
\begin{lem} \label{qnbeta0} (Kraaikamp and Liardet)
If there exists $n_0$ such that $\|q_n \beta\|\le {1\over 4} q_n \|q_n \alpha\|$ for $n \ge n_0$, then $\beta \in \Z \alpha + \Z$.
\end{lem}

The ratio
\begin{eqnarray}
c_n(\beta) := {\|q_n \beta\| \over q_n \|q_n \alpha\|} \label{defcn}
\end{eqnarray}
will be important in the proof of Theorem \ref{thmfunctEqu1}.
The previous lemma implies that, if $c_n(\beta) \leq \frac14, \forall n \geq n_0$ for some $n_0 \geq 1$, then $\beta \in \Z \alpha + \Z$.

We will use also the following lemma:
\begin{lem} \label{half} 1) If $q_n$ is odd and $q_n \|q_n\alpha\| < 1/2$, then $F^{(q_n)} = \pm 1$.

2) If $q_n$ is even and $q_n \|q_n\alpha\| < 1/2$, then $F^{(q_n)}(x) = \pm 2$ on a set $I_n$ of measure $\mu(I_n) \leq {1\over a_{n+1}}$
and $=0$ elsewhere.
\end{lem}
\proof The discontinuities of $F^{(q_n)}(x) = \sum_{j=0}^{q_n-1}F(x +j\alpha)$ are $t - j \alpha \mod 1$, with $j= 0,..., q_n-1$, $t = 0, \frac12$ and
the respective jumps are $+2, -2$.

1) Let us consider the case 1) where $q_n$ is odd. The discontinuities are of the form
${r\over q_n} - j_1(r) {\theta_n \over q_n}$, ${r\over q_n} - j_2(r) {\theta_n \over q_n} + {1\over 2q_n}$, with jumps $\pm 2$, where $0 \leq j_1(r), j_2(r) < q_n$,
for $r= 0, ... q_n-1$.
They belong respectively to $[{r \over q_n} - {\delta_n \over q_n}, \, {r \over q_n} + {\delta_n \over q_n}]$ and
$[{r \over q_n} + {1 \over 2 q_n} - {\delta_n \over q_n}, \, {r \over q_n} + {1 \over 2 q_n} + {\delta_n \over q_n}]$,
where $\delta_n := q_n\|q_n\alpha\| < 1/2$.

As ${r \over q_n} + {\delta_n \over q_n} < {r \over q_n} + {1 \over 2 q_n} - {\delta_n \over q_n}, $, the successive jumps of $F^{(q_n)}$ are alternatively
$+2, -2$, so that the values of $F^{(q_n)}$ are $u$ or $u - 2$, for a constant $u$. But $F$ is antisymmetric: $F(x+1/2)=-F(x)$, so also $F^{(q_n)}$ is antisymmetric
and non constant. Therefore the set of values is $\{u, u-2\} =\{-u, -u +2\}$, which implies $u = 1$.

2) Suppose $q_n$ even. For $r= 0, ..., q_n-1$, there are now two discontinuities (with jump respectively $+2$, $-2$ in any order)
in $[{r \over q_n} - {\delta_n \over q_n}, \, {r \over q_n} + {\delta_n \over q_n}]$. It shows that the set of values of $F^{(q_n)}$ belongs to
$\{u, u+2, u-2\}$, for some integer constant $u$.

There is a constant $u \in \Z$ such that, for each $r = 0, ..., q_n-1$, in restriction to the interval $[{r \over q_n} - {1 \over 2q_n}, \, {r \over q_n} + {1 \over 2q_n}]$,
the function $F^{(q_n)}$ takes the value $u + 2$ or $u - 2$ on a subinterval $I_{n, r}$ of length $\leq {|j_1(r) - j_2(r)| \over q_n} |\theta_n| \leq |\theta_n|$
and the value $u$ elsewhere.

Therefore, we have $0 = \int_{\T} F^{(q_n)} \, d\mu = u + 2 \sum_r \pm \mu(I_{n, r}$), hence $|u| \leq  2 \sum_r \mu(I_{n, r}) \leq 2 q_n |\theta_n| = 2 \delta_n < 1$.
This implies $u=0$. Let $I_n = \cup_r I_{n, r}$. As $\sum_r \mu(I_{n, r}) \leq q_n |\theta_n| \leq {1 \over a_{n+1}}$, the point 2) of the lemma is proved.
\eop

\section{\bf Proof of the first part of Theorem~\ref{Mainthm}}\label{sec1}

{\bf Proof of 1a)} \ The fact that if $\alpha$ has bounded partial quotients, then $T^f$ has (finite) Ratner's property has been proved in \cite{Fr-Le-Le}. Moreover, in \cite{Ka-Le} it has been proved that each flow satisfying (finite) Ratner's property has at most countable essential centralizer. On the other hand, we have already noticed that rigid flows have uncountable essential centralizer, whence our $T^f$ cannot be rigid.

{\bf Proof of 1b)} \  We want to show that
if $\alpha$ has bounded partial quotients, then, for a non-trivial $\beta$, equation~\eqref{funeq0} has no measurable solution. Our claim follows from the following result:

\begin{proposition} \label{zergodic} If $\alpha$ has bounded partial quotients, then for $\beta \not \in \Z \alpha + \Z$ the cocycle $\Phi_\beta$ is ergodic (as a $\Z$-valued cocycle).
\end{proposition}
\proof We use Lemma 2.3 and Proposition 3.8 in \cite{CoPi14}. The lemma shows that the cocycle $\Phi_\beta$ has ``well separated discontinuities'',
which implies ergodicity by the proposition. \eop

\section{\bf Rigidity of a special flow, proof of parts 2a), 2b) of Theorem~\ref{Mainthm}}\label{sec1a}

Assume that $(T_t)_{t\in \R}$ is a (measurable) measure-preserving flow on a probability standard Borel space $\zdr$.
We will consider $T_t$ as a Markov operator\footnote{Recall that  a linear contraction $\Phi$ on $L^2\zdr$
is called Markov, if $\Phi 1=\Phi^\ast 1=1$ and $\Phi h\geq0$ whenever $h\geq0$. The set of Markov operators is a convex set which is closed
(hence compact) in the weak operator topology.} on $L^2\zdr$: $T_tf:=f\circ T_t$. The following result is essentially due to V.~Ryzhikov (private communication).

\begin{Lemma}\label{l1} Assume that $(r_n)$ is a sequence of real numbers tending to $\infty$. Assume moreover that
$$T_{mr_n}\tend{n}{\infty} \frac12(T_{-m}+T_m)\text{ for all }m\in\Z,$$
weakly in the set of Markov operators. Then, the flow $(T_t)_{t\in\R}$ is rigid.
\end{Lemma}
\proof \  By assumption, for all $m\geq1$, we have
\beq\label{e1}
T_{m([r_n]+\{r_n\})}\tend{n}{\infty} \frac12(T_{-m}+T_m).
\eeq
By passing to a subsequence, if necessary, we have $T_{\{r_n\}}\tend{n}{\infty} T_r$ for some $r\in[0,1[$, whence
\beq\label{e2}
T_{m\{r_n\}}\tend{n}{\infty} T_{mr}\text{ strongly, for all }m\in\Z.\eeq
Since the convergence in \eqref{e2} is strong, by \eqref{e1}, we have
$$T_{m[r_n]}\tend{n}{\infty}\frac12(T_{-m}+T_m)\circ T_{-mr}\text{ for all }m\in\Z.$$
Using basic properties of the weak operator topology,\footnote{If $d$ is a metric compatible with the weak topology, then $d(C_n\circ A,C_n\circ B)\tend{n}{\infty}
d(C\circ A,C\circ B)$ for any linear contractions $A,B$ and $C_n\tend{n}{\infty}C$. The subsequence $(r_{n_k})$ is selected inductively, at the induction step,
$r_{n_{k+1}}$ is chosen so that $T_{j[r_{n_{k+1}}]}$ is so close to $D_j:=\frac12(T_{-j}+T_j)\circ T_{-jr}$ for $j=1,\ldots,k$
to have $T_{j([r_{n_{k+1}}]-[r_{n_k}])}$ is almost as close to $D_j\circ D_j^\ast$ as $T_{j[r_{n_k}]}$ is close to $D_j$, also $T_{(k+1)[r_{n_{k+1}}]}$
is very close to $D_{k+1}$.} we can now choose a sparse subsequence $(r_{n_k})$ so that, for all $m\in\Z$,
$$T_{m([r_{n_{k+1}}]-[r_{n_k}])}\tend{k}{\infty} \frac12(T_{-m}+T_m)\circ T_{-mr}\circ\left(\frac12(T_{-m}+T_m)\circ T_{-mr}\right)^\ast$$
$$=\frac12(T_{-m}+T_m)\circ \frac12(T_{-m}+T_m)=\frac14T_{-2m}+\frac12 Id+\frac14 T_{2m}.$$

By passing to a further subsequence, if necessary, we can assume that $[r_{n_k}]$ were chosen so that either they are all even or they are all odd.
This yields
\beq\label{e3}
T_{2m j_k}\tend{k}{\infty} \frac14T_{-2m}+\frac12 Id+\frac14 T_{2m}, \text{ for all }m\in\Z,
\eeq
where  $j_k=([r_{n_{k+1}}]-[r_{n_k}])/2$.

It follows that for each $\ell\in \Z$ the operator $\frac14T_{-2\ell}+\frac12 Id+\frac14 T_{2\ell}$ is an accumulation point  of the set $\{T_{n}:\:n\in\Z\}$.
Fix $m\geq1$ and $k\geq1$. By taking $\ell=mj_k$, it follows  that the operator $\frac14T_{-2mj_k}+\frac12 Id+\frac14 T_{2mj_k}$ is an accumulation point
of $\{T_{n}:\:n\in\Z\}$. Letting $k\to\infty$ and using~\eqref{e3}, we obtain that the operator
$$\frac14\left(\frac14 T_{-2m}+\frac12 Id+\frac14 T_{2m}\right)+\frac12 Id+ \frac14\left(\frac14T_{2m}+\frac12 Id+\frac14 T_{-2m}\right)=
\frac18T_{2m}+\frac34 Id+\frac18 T_{-2m}$$
is an accumulation point of the set $\{T_n:\:n\in\Z\}$. By iterating this procedure, we obtain that $Id$ is an accumulation point of $\{T_n:\:n\in\Z\}$ and the result follows.
\eop

Consider now $R_\alpha x=x+\alpha \mod 1$ an irrational rotation on $\T$. Recall that $(q_n)$ denote the sequence of denominators of $\alpha$. Recall also that,
 for any function $\varphi$ on $\T$ and a positive integer $\ell$,  we denote by $\varphi^{(\ell)}$ the ergodic sum
$\varphi^{(\ell)} (x) = \sum_{k=0}^{\ell - 1}\varphi(x + k \alpha)$ (cf.\ \eqref{sumy}).

Let $f:\T\to\R^+$ be of bounded variation. As noticed in \cite{Le-Pa}, we have
\beq\label{e4}
\|f^{(mq)}-mf_q\|_\infty \leq \frac12 m^2q\|q\alpha\|{\rm Var}(f),
\eeq
where $q$ is a denominator of $\alpha$ and $f_{\ell}$ is the periodized function $f_{\ell}(x)=\sum_{i=0}^{\ell-1} f(x+{i \over \ell})$.

Assume that $\alpha$ has unbounded partial quotients and let
$q_{n_k}\|q_{n_k}\alpha\|\to 0$ along some subsequence $(q_{n_k})$ of the sequence $(q_n)$ of denominators of $\alpha$. Set $ c:=\int_X f \, d\mu$ and $F:=f - c$.

We can assume additionally that $(F^{(q_{n_k})})_\ast\tend{k}{\infty} P$ in distribution ($P$ is a probability measure concentrated on $[-{\rm Var}(f),{\rm Var}(f)]$
by the Denjoy-Koksma inequality\footnote{Recall that the Denjoy-Koksma inequality states
$|f^{(q_n)}(x)|\leq {\rm Var}(f)$ for each zero mean, bounded variation $f:\T\to\R$, $n\geq1$ and $x\in\T$.}). Denoting by $mP$ the image of $P$ via the map $r\mapsto mr$, it follows by~\eqref{e4} that
\beq\label{e5} (F^{(mq_{n_k})})_\ast \tend{k}{\infty} mP \eeq
in distribution for each $m\in \Z$. By \cite{Fr-Le}, we hence obtain the following weak convergence in the space of Markov operators:
\beq\label{e6}
T^f_{mcq_{n_k}}\tend{k}{\infty} \int_{\R}T^f_{-t}\,d(mP)(t).\eeq

Consider now our special case (cf.\ \eqref{ab}) of $f=f_{a,b}$ for which
$f(x)=a$ for $x\in [0,\frac12)$ and $f(x)=b$ for $x\in[\frac12,1[$. We assume that $a,b>0$. Then $c=\frac12(a+b)$ and,
if moreover we take $a-b = 2$, $F$ now becomes
$$F= 1_{[0, \frac12[} - 1_{[\frac12, 1[}.$$

Using Lemma \ref{half} (Section~\ref{Prelim}), the following immediately follows:
\begin{Lemma}\label{l3} Assume moreover that the denominators $q_{n_k}$ above are all odd and $a - b = 2$. Then
$(F^{(q_{n_k})})_\ast \   {\underset{k \to \infty} \longrightarrow}  \ \frac12(\delta_{-1}+\delta_1)$.
\end{Lemma}
It follows from \eqref{e5} that
$\displaystyle (F^{(mq_{n_k})})_\ast  \tend{k}{\infty} \ \frac12(\delta_{-m}+\delta_m)$
in distribution for each $m\in\Z$ and then by \eqref{e6} that
$$T^f_{mcq_{n_k}} \  {\underset{k \to \infty} \longrightarrow}  \ \frac12(T^f_{-m}+T^f_m)$$
weakly in the set of Markov operators, for each $m\in\Z$.

\vskip 3mm
{\bf Proof of part 2a) and 2b) of Theorem~\ref{Mainthm}}

Assume that $\alpha$ has unbounded partial quotients.

Then, either there is a subsequence $(n_k)$ such that  $q_{n_k}\|q_{n_k}\alpha\|\tend{k}{\infty}0$,
where each denominator $q_{n_k}$, $k\geq1$, is odd. If $a - b = 2$, then the special flow $T^f$, obtained by $T = R_\alpha$ and $f=f_{a,b}$, is rigid.
Indeed, the result follows from the previous discussion, using Lemma~\ref{l1} (with $r_n=cq_n$).

Or, there is a subsequence $(n_k)$ such that  $q_{n_k}\|q_{n_k}\alpha\|\tend{k}{\infty}0$, where each denominator $q_{n_k}$, $k\geq1$, is even.
Then, by the second part of Lemma~\ref{half}, it implies that $F^{(q_{n_k})}\to 0$ in $L^2$. The result follows by a folklore argument
(cf. Proposition \ref{series1} in Appendix and the remark below).

This shows 2a). Part 2b) follows also from what precedes. \eop

\vskip 3mm
{\bf Remark on rigidity}

Let us consider a special flow over the rotation $R_\alpha$ by an irrational $\alpha$, under a roof function $f$ in $L^2$. Assume that $\int f \, d\mu = 1$.
Let  $f_0$ denote the centered function $f - 1$. It is not hard to see that the existence of a sequence $(r_n)$ of integers tending to infinity such that
\begin{eqnarray}
\ \|r_n \alpha\| \to 0, \ \|f_0^{(r_n)}\|_2 \to 0 \label{rigid1}
\end{eqnarray}
implies the rigidity of the special flow $(T^f)$ (cf. proof of Proposition~\ref{pot4} in Appendix).

The following lemma shows that for $f_0 = \varphi := 1_{[0, {1 \over 2}[} - 1_{[{1 \over 2},1[}$ there is no sequence $(r_n)$ satisfying~(\ref{rigid1})
if $q_k$ is odd for $k$ big enough. This implies that that the method given by~(\ref{rigid1}) cannot be used to prove rigidity in the framework
of Theorem~\ref{Mainthm}, case 2b).

\begin{lem}\label{exclusive} There is $\delta > 0$ such that, if $q_k$ is odd for $k \geq k_0$, then for every integer $s \geq k_0$, $\|\varphi^{(s)}\|_2\geq \delta$.
\end{lem}
\proof \  We have $\varphi(x)  = \sum_{r\in\Z} {2 \over \pi i (2r+1)} \, e^{2\pi i (2r+1) x}$, hence (cf.\ \eqref{es2} below)
\begin{eqnarray*}
&&\|\varphi^{(s)}\|_2^2 =  {4 \over \pi^2} \, \sum_{r\in\Z} \, {1 \over (2r+1)^2} \, \left({\sin \pi s (2r+1) \alpha \over \sin \pi (2r+1) \alpha}\right)^2.
\end{eqnarray*}
There is $k$ such that $q_k \leq s < q_{k+1}$. Taking the term corresponding to $2r+1 = q_k$, in the above series and using the equality
$\|sq_k\alpha\| = s\|q_k\alpha\|$ valid since $s < q_{k+1} \leq 1/ \|q_k \alpha\|$, we get (up to a constant):
\begin{eqnarray*}
&&\|\varphi^{(s)}\|_2 \geq {1\over q_k} \, {s \|q_k \, \alpha\| \over \|q_k \alpha\|} = {s \over q_k} \geq1.
\end{eqnarray*}
\eop

\section{\bf Centralizer and functional equation for $\Phi_\beta$}\label{sec2}

\subsection{\bf Regularity of a class of step cocycles}

\

In this subsection it is shown the existence of a large class of values of $\beta$ such that Property~(II) defined in Section \ref{Prelim} holds for $\Phi_\beta$.
We start by a general result based on a Fourier computation.

If $\varphi$ is a centered BV function, we write $\varphi (x) = \sum_{r \not = 0} {\gamma_r(\varphi) \over r} \, e^{2 \pi i r x}$ for its Fourier series
and we have $\sup_r |\gamma_r(\varphi)| < \infty$.
\begin{thm} \label{FourMeth} Let $\varphi$ be a centered BV real valued function. If there are a subsequence $(q_{n_k})$ of denominators
and a constant $\delta > 0$ such that
\begin{eqnarray}
&&|\gamma_{q_{n_k}}(\varphi)| \geq \delta, \forall k \geq 1, \label{condgammak} \\
&&M:=\sup_{k \,: \, a_{n_k+1}=1} a_{n_k} < \infty, \label{condak}
\end{eqnarray}
then the cocycle generated by $\varphi$ has a finite essential value $\not = 0$ (hence  $\varphi$ is regular and is not a coboundary).
\end{thm}
\proof We claim that there is a positive constant $c$ such that $\|\varphi^{(q_{n_k})}\|_2^2 \geq c, \, \forall k \geq 1$.

By Lemma \ref{supp}, since by the Denjoy-Koksma inequality $\varphi^{(q_n)}$ is uniformly bounded by $\Var(\varphi)$,
this will imply that the cocycle generated by $\varphi$ has a non zero essential value,
hence is regular. Moreover, since $\cal E(\varphi) \not = \{0\}$, $\varphi$ is not a coboundary.

Now, we prove the claim. The ergodic sum of $\varphi$ at time $q$ and the square of its $L^2$-norm read:
\begin{eqnarray}
&&\varphi^{(q)}(x) = \sum_{r \not = 0} {\gamma_{r}(\varphi) \over r}  \, { e^{2\pi i q r \alpha} - 1 \over e^{2 \pi i r \alpha} - 1}
\, e^{2\pi i r x}, \label{es1} \\
&&\|\varphi^{(q)}\|_2^2 =  \sum_{r \not = 0} {|\gamma_{r}(\varphi)|^2 \over r^2} \, \left({\sin \pi q r \alpha \over \sin \pi r \alpha}\right)^2. \label{es2}
\end{eqnarray}
Taking the term corresponding to $r = q_n$, we have, for $n$ in the sequence $\Cal S = (n_k)$:
\begin{eqnarray} \label{e:jp1}
&&\|\varphi^{(q_n)}\|_2^2 \geq \delta^2 \,  {1\over q_n^2} \, \left({\sin \pi q_n^2 \, \alpha \over \sin \pi q_n \alpha}\right)^2
= \delta^2 \, {1\over q_n^2} \, \left({\sin \pi q_n \, \theta_n \over \sin \pi \theta_n}\right)^2,
\end{eqnarray}
with (see \eqref{f_3b})
$$q_n \, |\theta_n| \leq {q_n \over a_{n+1} q_n + q_{n-1}} \leq  {1 \over a_{n+1}}.$$
If $a_{n+1} \geq 2$, then $q_n \, |\theta_n| \leq \frac12$ and it follows from (\ref{sinA}) and~\eqref{e:jp1} that:
\begin{eqnarray}
&&\|\varphi^{(q_n)}\|_2^2 \geq \delta^2 \, {1\over q_n^2} \,
\left({2 q_n \, \theta_n \over \pi \theta_n}\right)^2 = 4 {\delta^2 \over \pi ^2}.
\end{eqnarray}
Now, for $n$ in $\Cal S = (n_k)$, suppose that $a_{n+1} = 1$, so that $q_{n+1} = q_n + q_{n-1}$. Considering still $r= q_{n}$,
but in the Fourier series of $\varphi^{(q_{n-1})}$, we get the lower bound:
\begin{eqnarray}
&&\|\varphi^{(q_{n-1})}\|_2^2 \geq \delta^2 \, {1\over q_{n}^2} \,
\left({\sin \pi q_{n-1} \, q_n \alpha\over \sin \pi q_{n} \alpha}\right)^2 = \delta^2 \, {1\over q_{n}^2} \,
\left({\sin \pi q_{n-1} \, \theta_{n} \over \sin \pi \theta_{n}}\right)^2.
\end{eqnarray}
By \eqref{f_3} we have $q_{n-1} \, |\theta_{n}| \leq {q_{n-1} \over q_{n+1}} = {q_{n-1} \over q_{n} + q_{n-1}}\leq \frac12$,
so we can use (\ref{sinA}). From the hypothesis, for this value of $n$, we have $a_n \leq M$, so that
$$ {1 \over q_n} \, {|\sin \pi q_{n-1} \, \theta_{n}| \over |\sin \pi \theta_{n}|}
\geq {2\over \pi} {q_{n-1} \over q_{n}} \geq {2\over \pi} {q_{n-1} \over a_n q_{n-1} + q_{n-2}} = {2\over \pi}{1 \over a_n + q_{n-2}/q_{n-1}}
\geq {2\over \pi} \, {1 \over M + 1}. \eop$$

{\it Remark}:  When $\varphi$ has values in $\Z$ as in the examples below, we can give the following variant of the previous proof.
By the Denjoy-Koksma inequality, $\varphi^{(q_n)}$ takes a finite number of integral values in $[-\Var(\varphi), \Var(\varphi)]$ and
$$\|\varphi^{(q_n)}\|_2^2 = \sum_{j: |j| \leq \Var(\varphi)} j^2 \mu( \{\varphi^{(q_n)} = j\}).$$

By the claim in the previous proof, $\|\varphi^{(q_{n_k})}\|_2^2 \geq c$ for a positive constant $c$. This implies that, on a subsequence of $(n_k)$,
$\varphi^{(q_{n_k})}$ takes a fixed value $j_0 \not = 0$ on sets whose measure is bounded away from 0.
Therefore, $j_0$ is a non zero essential value.

\vskip 3mm
\subsubsection{\bf Examples, the step function $\Phi_\beta$}

\

Now, we consider specific examples related to Theorem \ref{Mainthm} and introduce some notation.
The argument of the functions below are understood $\mod 1$.

Let $G(x) = \{x\} -\frac12$, $F = 1_{[0, \frac12[} - 1_{[\frac12, 0[}$, as above, and $\Phi_\beta := \frac12 F(.-\beta) - \frac12 F$,
where $\beta$ is a real number in $]0, 1[$.

If $0 < \beta \leq \frac12$, then $\Phi_\beta = -1_{[0, \beta[} + 1_{[\frac12, \beta + \frac12[}$;
if $\frac12 < \beta <1$, then $\Phi_\beta = - 1_{[\beta - \frac12, \frac12[} + 1_{[\beta, 1[}$.

For $\beta \not = \frac12$, the jumps of $\Phi_\beta$ are respectively $+1$ at $\beta$, $+1$ at $\frac12$, $-1$ at $\frac12 + \beta \mod 1$.
The jump at 0 is $\lim_{t \to 0^+} \Phi_\beta(t) - \lim_{t \to 1^-} \Phi_\beta(t) = -1$.

Observe also that $F= 2(R_\frac12 -I) G$, so $\Phi_\beta : = \frac12(R_{-\beta} - I) F = (R_{-\beta} - I) (R_{-\frac12} - I)G$.

More generally, let $\beta_1, ..., \beta_v$ be real numbers and set
\begin{eqnarray}
&&\varphi_{\beta_1, ..., \beta_v} := \prod_{j= 1}^v \, (R_{-\beta_j} - I) \ G. \label{phibetaj}
\end{eqnarray}
With this notation, the function $\Phi_\beta$ considered before is $\varphi_{\frac12, \beta}$ and we have
$\varphi_{\beta, \gamma} = - (I - R_{-\gamma}) \, \zeta_\beta$, with $\zeta_\beta := 1_{[0, \beta[} - \beta$ since $(R_{-\beta}-I)G=\zeta_\beta$.
The Fourier series of $G$ and $\varphi_{\beta_1, ..., \beta_v}$ are respectively
\begin{eqnarray*}
G(x) = {-1 \over 2\pi i} \, \sum_{r \not = 0} {1\over r} \  e^{2\pi i r x},
\ \ \varphi_{\beta_1, ..., \beta_v}(x) = {-1 \over 2\pi i} \, \sum_{r \not = 0} {1\over r} \ \prod_{j= 1}^v \, (e^{-2\pi i r \beta_j} - 1)
\, e^{2\pi i r x},
\end{eqnarray*}
and therefore
$$|\gamma_{q_n}(\varphi_{\beta_1, ..., \beta_v})| = {2^{v-1} \over \pi} \, \prod_{j= 1}^v \, |\sin \pi q_n \beta_j|.$$

Immediately from Theorem \ref{FourMeth}, we obtain the following results:
\begin{cor} \label{corphibeta0} If $\limsup_n \prod_{j=1}^v \, \|q_n \beta_j\| >0$ and $\sup_{n: a_{n+1}=1} a_n < \infty$, then the group of finite
essential values of the cocycle $\varphi_{\beta_1, ..., \beta_v}$ is not reduced to 0.
\end{cor}

\begin{cor} \label{corphibeta}  If there is subsequence $(q_{n_k})_{k \geq 1}$ such that  $q_{n_k}$ is odd,  $\sup_{k \,: \, a_{n_k+1}=1} a_{n_k} < \infty$
and $\limsup_k \|q_{n_k} \beta\| > 0$, then $\Phi_\beta = \varphi_{\frac12, \beta}$ is regular and is not a coboundary.

In particular, this is true if $q_n$ is odd for $n$ big enough and $\limsup_n \|q_n \beta\| > 0$.
\end{cor}
\proof The particular case follows from Remark \ref{oddeven1}, which shows that, if  $q_n$ is odd for $n \geq n_0$, then $a_n$ is even, hence $\geq 2$,
for $n \geq n_0 +1$. \eop

\vskip 3mm
\section{\bf Coboundary equation, end of the proof of Theorem~\ref{Mainthm}}\label{sec3}

The aim of this section is to finish the proof of Theorem~\ref{Mainthm} by proving 2c) and 2d). We start by a preliminary discussion on the discontinuities of
$\varphi_\beta^{(q_n)}$ for a general $\varphi$ and then of $\Phi_\beta^{(q_n)}$, which will be useful in the proof of Theorem \ref{thmfunctEqu1}.

Let $\gamma$ be in $[0, 1[$ and $n\geq 1$. Recall that $\|q_n\gamma\| \leq \frac12$. We define $t(\gamma, n) \in \Z$ and $\varepsilon_n(\gamma) = \pm 1$ by
\begin{equation}
q_n \gamma = t(\gamma, n) + \varepsilon_n(\gamma) \|q_n\gamma\|.  \label{qngamma}
\end{equation}

So, we have:
\begin{eqnarray}
\gamma = {t(\gamma, n) \over q_n} + \varepsilon_n(\gamma) {\|q_n\gamma\| \over q_n}.\label{qngamma3}
\end{eqnarray}
Note that if $\gamma = \alpha$, then (\ref{qngamma3}) reads (cf.\ the last equality in
(\ref{f_5})):
\begin{eqnarray}
\alpha = {p_n \over q_n} + (-1)^{n} {\|q_n \alpha\| \over q_n}. \label{f_5z}
\end{eqnarray}
If $\gamma = \frac12$, then $\|q_n \frac12\| = 0$ or $\frac12$, depending whether $q_n$ is even or odd.

Suppose now that $q_n$ is odd: $q_n = 2 q_n' +1$. Then, by~\eqref{qngamma3}, since $\frac12 = {q_n' \over q_n} + {1 \over 2 q_n}$, we have:
\begin{equation}\label{aa2}
\gamma + \frac12 ={t(\gamma, n) + q_n' \over q_n} + {1 \over 2 q_n} + \varepsilon_n(\gamma) {\|q_n\gamma\| \over q_n}.
\end{equation}

In (\ref{qngamma}),  $t(\gamma, n)$ and $\varepsilon_n(\gamma)$  are uniquely defined, excepted for $\gamma = \frac12$. For this special value,
$\|q_n \frac12\| = \frac12$ and we have the representation $\frac12 = {q_n' \over q_n} +{ {1 \over 2 } \over q_n}$.

\vskip 3mm
{\it Location of the discontinuities of $\varphi^{(q_n)}$}

Let $\varphi$ be a 1-periodic function. If $\gamma$ is a discontinuity of $\varphi$,
the discontinuities of $\varphi^{(q_n)}$corresponding to $\gamma$ are located at $\gamma - \ell \alpha \mod 1$, $\ell = 0, 1, ..., q_n -1$.
We call them {\em discontinuities of type} $\gamma$.

For a given denominator $q_n$, we consider the grid $\{0, {1\over q_n}, {2\over q_n}, ..., {q_n-1 \over q_n}\}$ and denote by
$I_{n,k} = I_k$ the interval $[{k\over q_n}, {k+1\over q_n}[$, $0 \leq k < q_n$. In each interval $I_k$, there is one and only
one discontinuity of type 0. For $0 < \gamma < 1$, there are 0, 1 or 2 discontinuities of type $\gamma$ in each interval $I_k$,
since $\|\ell \alpha\| > \|q_n \alpha\| \geq {1 \over 2 q_n}$ for $\ell = 1, ..., q_{n}-1$,  by \eqref{f_4}.

For $n \geq 1$ and $\gamma \in [0, 1[$, the map $\ell \to - \ell p_n + t(\gamma, n) \mod q_n$ defines a permutation
of the set $\{0, 1, ..., q_n-1\}$. In view of~\eqref{aa1}, its inverse map is $k \to u_n(k, \gamma)$, where
\begin{equation}
u_n(k, \gamma) \in \{0, 1, ..., q_n-1\}  \text{ and } u_n(k, \gamma) = (-1)^{n-1} q_{n-1} (- k + t(\gamma, n)) \mod q_n. \label{unkdef}
\end{equation}
We put $A_n(\gamma) = \varepsilon_n(\gamma) \|q_n \gamma\|$, $B_n(\gamma, k) = (-1)^{n-1} u_n(k, \gamma) \, \|q_n \alpha\|$.

Using \eqref{qngamma3}, \eqref{f_5z} and the definition of $u_n(k,\gamma)$, for each discontinuity $\gamma$ of $\varphi$,
we can label the discontinuities $\gamma - \ell \alpha \mod 1$, $\ell = 0, ..., q_n-1$, as $\zeta(\gamma, k, n)$, $k = 0, ..., q_n-1$:
\begin{eqnarray}
\zeta(\gamma, k, n) &=&{k \over q_n} + {\varepsilon_n (\gamma) \|q_n\gamma\| \over q_n}  + { (-1)^{n-1} u_n(k, \gamma)
\, \|q_n\alpha\| \over q_n} \label{locak1}\\
&=& {k \over q_n} + {A_n(\gamma) + B_n(\gamma, k) \over q_n}. \label{locak2}
\end{eqnarray}
We have $\|q_n \gamma\| \leq \frac12$, $ u_n(k, \gamma) \, \|q_n \alpha\| < q_n \|q_n \alpha\| < 1/a_{n+1} \leq \frac12$, since $a_{n+1}$ is even
(because the $q_n$'s are odd). Hence $|A_n(\gamma) + B_n(\gamma, k)| < 1$.

Thus, equation (\ref{locak1}) gives the position of the discontinuities of type $\gamma$:
$\zeta(\gamma, k, n)$ belongs to the interval $I_k$ if $A_n(\gamma) + B_n(\gamma, k) >0$, to the interval $I_{k-1}$ if $A_n(\gamma) + B_n(\gamma, k) <0$.
(with the convention $I_{-1} = I_{q_{n -1}}$). Moreover,
for each $\gamma$, the sequence $(\zeta(\gamma, k, n)), k= 0, ..., q_n-1)$ is increasing, since:
$$\zeta(\gamma, k+1, n) - \zeta(\gamma, k, n) = {1\over q_n} +(-1)^{n-1} {\|q_n \alpha\| \over q_n} \, (u_n(k+1, \gamma) - u_n(k, \gamma))$$
and $\|q_n \alpha\| \, |u_n(k+1, \gamma) - u_n(k, \gamma)| \leq q_n \|q_n \alpha\| < 1$.

\vskip 3mm
{\it Discontinuities of $\Phi_\beta^{(q_n)}$}

The discontinuities of $\Phi_\beta^{(q_n)}$ are of type $0$, $\frac12$, $\beta$ and $\beta + \frac12$.
For the type $0$, $\frac12$ and $\beta$, they read, for $k= 0, 1, ..., q_n-1$,
\begin{eqnarray}
\zeta(0, k, n) &&= {k \over q_n} + (-1)^{n-1} {u_n(k, 0)  \, \|q_n\alpha\| \over q_n}, \label{disczer}\\
\zeta(\frac12, k, n) &&=  {k \over q_n} + {1 \over 2q_n} + (-1)^{n-1} {u_n(k, \frac12)  \, \|q_n\alpha\| \over q_n}, \label{dischalf}\\
\zeta(\beta, k, n) &&=  {k \over q_n} + \varepsilon_n(\beta){\|q_n\beta\| \over q_n} + (-1)^{n-1} {u_n(k, \beta)  \, \|q_n\alpha\| \over q_n}.
\label{discbet1}
\end{eqnarray}

The discontinuities of type $\beta + \frac12$ can be written
\begin{eqnarray}
\zeta_1( \beta + \frac12, k, n) &&= {k \over q_n} + {1 \over 2q_n} +  \varepsilon_n(\beta) {\|q_n \beta\|\over q_n}
+ (-1)^{n-1} {u_n'(k, \beta)  \, \|q_n\alpha\| \over q_n}, \label{discbeta2}
\end{eqnarray}
where $u_n'(k, \beta) = u_n(k, \beta) + (-1)^{n-1} q_{n-1} q_n' \mod q_n$. Indeed, using \eqref{discbet1} and  $\frac12 = {q_n' \over q_n} + {1 \over 2 q_n}$,
we have (mod 1)
$$\zeta(\beta, \ell, n) + \frac12 =  {\ell + q_n'  \over q_n} + {1 \over 2 q_n} + \varepsilon_n(\beta){\|q_n\beta\| \over q_n} + (-1)^{n-1} {u_n(\ell, \beta)
\, \|q_n\alpha\| \over q_n}$$
and by taking $\ell = k -q_n'$, we get  \eqref{discbeta2}.

Let us assume $n$ {\it odd} (hence $\alpha = {p_n \over q_n} - {\|q_n\alpha\| \over q_n}$). The discussion is analogous for $n$ even.

For $\gamma = 0$, there is one and only one discontinuity, namely $\zeta(0,k,n)$, of type 0 of $\Phi^{(q_n)}$ in $I_k$. By \eqref{disczer}, it reads
$\displaystyle {k \over q_n} + { u_n(k, 0)  \, \|q_n\alpha\| \over q_n}$. Since $u_n(k,0) \, \|q_n\alpha\| \leq q_n \, \|q_n\alpha\| \leq {1 \over a_{n+1}}$,
this discontinuity is close to the left endpoint ${k \over q_n}$ of $I_k$ if $a_{n+1}$ is big.

By (\ref{dischalf}), the discontinuity of type $\frac12$ in $I_k$ is $\displaystyle \zeta(\frac12, k,n)= {k \over q_n} + {1 \over 2q_n} + {u_n(k, \frac12)
\, \|q_n\alpha\| \over q_n}$, with $0 \leq u_n(k, \frac12)  \, \|q_n\alpha\| \leq q_n \, \|q_n \alpha\| < {1 \over a_{n+1}}$, hence located close
to the middle of $I_k$. It is the only discontinuity of type $\frac12$ belonging to $I_k$ if $a_{n+1}$ is big.

By (\ref{discbet1}), $\zeta(\beta, k, n)$, discontinuity of type $\beta$, is close to ${k \over q_n}$ (hence close to $\zeta(0, k, n)$, either to the left
or to the right of it), if $\|q_n\beta\|$ is small and $a_{n+1}$ is big. Furthermore, notice that the next discontinuity of type $\beta$, $\zeta(\beta, k+1, n)$,
may belong to $I_k$, but is close to the right endpoint ${k+1 \over q_n}$ of $I_k$.

By (\ref{discbeta2}), $\zeta_1(\beta + \frac12, k, n)$, discontinuity of type $\beta + \frac12$, is close to ${k \over q_n} + {1 \over 2q_n}$
(hence to $\zeta(\frac12, k, n)$, left or right), if $\|q_n\beta\|$ is small and $a_{n+1}$ is big.

We conclude these preliminaries by the following remark:
\begin{rem} \label{Remcob1}
The set $\{\beta\in\T: \Phi_\beta \text{ is an } R_\alpha-\text{coboundary} \}$ is an additive group, cf.~\eqref{e:lift0}.
\end{rem}

\begin{thm} \label{thmfunctEqu1} Assume that the denominators $q_n$ of $\alpha$ are odd for $n\geq n_0$, for some $n_0$. Then,
if $\beta \not \in \Z \alpha + \Z$, $\Phi_\beta$ is not a coboundary, i.e.,\ the functional equation
\begin{eqnarray}
\Phi_\beta(x) = g(x+\alpha) - g(x), \text{ for }\mu-\text{a.e. } x\in\T,\label{funeq1}
\end{eqnarray}
has no measurable solution~$g$.
\end{thm}
\proof Let us assume $\beta \not \in \Z \alpha + \Z$. For $\beta = \frac12$,  we get $\Phi_\frac12 = - F(x)$ which is not a coboundary for $R_\alpha$.
So we can assume $\beta \not = \frac12$.

For $n \geq n_0$, since $q_n$ is odd, we have $q_n \|q_n \alpha\| < a_{n+1}^{-1} \leq \frac12$.
Therefore by Lemma \ref{half}, $\sum_{j=0}^{q_n-1}F(x +j\alpha)=\pm 1$ for all $x$. It follows:
\begin{eqnarray}
\Phi_\beta^{(q_n)}(x)  = \frac12\sum_{j=0}^{q_n-1}F(x  - \beta +j\alpha)
- \frac12 \sum_{j=0}^{q_n-1}F(x +j\alpha)= 1, -1 \text{ or } 0. \label{sum01}
\end{eqnarray}

There are two cases depending on the behaviour of $\|q_n \beta \|$:

{\bf A)} $\limsup_ n \|q_n \beta\| > 0$

In this case, we use Corollary \ref{corphibeta} to conclude that $\Phi_\beta$ is not a coboundary. A stronger conclusion is the following:
there is a sequence $(n_k)$ and $\delta > 0$ such that $\|\Phi_\beta^{(q_{n_k})}\|_2 \geq \delta$ which implies by (\ref{sum01})
(cf.~Theorem~\ref{FourMeth} and Corollary~\ref{corphibeta}) that $\Phi_\beta^{(q_{n_k})}(x) = \pm 1$ on sets whose  measure is bounded away from 0.
Hence~1 is an essential value of the cocycle and the skew map $R_{\alpha, \Phi_\beta}$ is ergodic on $X \times \Z$. A fortiori, $\Phi_\beta$ is not a coboundary.

\vskip 3mm
{\bf B)} $\|q_n \beta \| \to 0$

We are going to show that in case B), for $\beta \not \in \Z \alpha + \Z$, the cocycle $\Phi_\beta$ is not a coboundary,
which will conclude the proof of the theorem. But contrary to case A), ergodicity of the skew product may fail (see Remark \ref{nonreg1} below).

We start by studying the support of $\Phi_\beta^{(q_n)}$ deduced from the location of the discontinuities of the cocycle as studied above.

\vskip 3mm
{\it Clusters of discontinuities and support of $\Phi_\beta^{(q_n)}$}

Let $n$ be such that $\|q_n \beta\|$ is small and $a_{n+1}$ is big. Then the picture is the following.

In the interval $I_k = [{k\over q_n}, {k+1 \over q_n}[$, there is one and only one discontinuity of type 0, $\zeta(0, k, n)$, and one
and only one of type $\frac12$, $\zeta(\frac12, k, n)$, which are close respectively to the left endpoint ${k \over q_n}$
and the middle point ${k \over q_n}$ + ${1\over 2q_n}$ of $I_k$.

There is a discontinuity of type $\beta$, $\zeta(\beta, k, n)$, at left or at right of $\zeta(0, k, n)$ and close to it.
Likewise, there is a discontinuity $\zeta(\beta + \frac12, k, n)$ of type $\beta+\frac12$
at left or at right of $\zeta(\frac12, k, n)$ and close to it.

The discontinuity which is the nearest discontinuity to the discontinuity $\zeta(0, k, n)$ of type 0 in $I_k$ is $\zeta(\beta, k, n)$
of type $\beta$. The nearest discontinuity close to a discontinuity of type $\frac12$ in $I_k$ is $\zeta(\beta + \frac12, k, n)$ of type
$\beta + \frac12$.

This shows that the discontinuities of $\Phi_\beta^{(q_n)}$ gather in well separated ``clusters'' (which here are groups of two discontinuities
close together) (see Fig 1: graph of $\Phi_\beta$,  with $\alpha=\pi -3$, $\beta= 2 -\sqrt 2$ and Fig. 2:  graph of $\Phi_\beta^{(q_1)}$, with $q_1 = 7$,
first denominator in the sequence of denominators of $\alpha$).

This situation can be described as follows. Suppose for concreteness that $\zeta(\beta, k, n)$ is located at the right to $\zeta(0, k, n)$.
If a point $x$ moves in $I_k$ to the right, starting close to $\zeta(0, k, n)$ at its left (hence close to ${k \over q_n}$), it crosses successively
two discontinuities $\zeta(0, k, n)$, $\zeta(\beta, k, n)$, with jumps respectively $-1, +1$. The cocycle $\Phi_\beta^{(q_{n})}(x)$ has a constant
value $v$ at the left of $\zeta(0, k, n)$, then again $v$ after the discontinuities and keeps this value until it is close to
${k\over q_n} +{1 \over 2 q_n}$.
Therefore it keeps a constant value $v$ on an interval of length close to ${1 \over 2q_n}$.
Then, still with $x$ moving in $I_k$, $\Phi_\beta^{(q_n)}(x)$ takes again the value $v$ after crossing the discontinuities
$\zeta(\frac12, k, n)$ and $\zeta(\beta+\frac12, k, n)$. It keeps this constant value on an interval again of length close to ${1 \over 2q_n}$.
On the whole interval $[0, 1]$, $\Phi_\beta^{(q_n)}$ takes the values $v$ on a set of measure close to 1.
Since the integral is 0, this implies $v=0$.

By Lemma \ref{qnbeta0}, we know that there is a subsequence $(n_j)$ such that $c_{n_j}(\beta) > \frac14$, where $c_n$ is the ratio (\ref{defcn}).
Since $\|q_{n_j} \beta\| \to 0$, from this and \eqref{f_3b}, it follows:
$4\|q_{n_j} \beta\| \geq q_{n_j} \|q_{n_j} \alpha\| \geq q_{n_j} (q_{n_j +1} + q_{n_j})^{-1} \geq (a_{n_j +1} + 2)^{-1}$.
Therefore, $a_{n_j +1} \to \infty$ since $\|q_{n_j} \beta\| \to 0$.

Since $\|q_{n_j} (8\beta)\| =  8 \|q_{n_j} \beta\|$ (because $\|q_{n_j} \beta\|$ is close to 0), replacing $\beta$ by $\beta' = 8\beta$,
we get $c_{n_j}(\beta') = c_{n_j}(8\beta) \geq 2$. Observe that if (\ref{funeq1}) has  a measurable solution for $\beta$,
then the equation (\ref{funeq1}) corresponding to $\beta'$ has a measurable solution, by Remark \ref{Remcob1}.

Therefore, for the proof of the non-existence of a measurable solution of (\ref{funeq1}), we can assume that,
for a strictly increasing sequence $\cal S= (n_j)$, we have $\|q_{n_j}\beta\| \to 0$ and $c_{n_j}(\beta) \geq 2$.

On the first half of the interval $I_k$, $\Phi_\beta^{(q_{n_j})}$ has its support on a small interval since $\|q_{n_j} \beta\|$ is small.
The idea of the proof is to consider $\Phi_\beta^{(L_j q_{n_j})}$ the cocycle at time $L_j q_{n_j}$ for a well chosen integer $L_j$.
(see Fig. 3: graph of $\Phi_\beta^{(3q_1)}$).

The idea of the proof is as follows. Suppose for concreteness that $n_j$ is odd. The support of $\Phi_\beta^{(L_j q_{n_j})}$ on the first half of $I_k$ is a union
of translates by multiples of $\theta_{n_j} = - \|q_{n_j} \alpha\|$ of the support of $\Phi_\beta^{(q_{n_j})}$. Up to a certain amount of translates,
there is no interference with the part of the support where $\Phi_\beta^{(q_{n_j})}$ has an opposite sign. To cover a set of measure $\geq \delta /q_{n_j}$
for some $\delta > 0$ inside the interval $I_{n_j}$, we take $L_j \sim \delta a_{n_j +1}$. This is enough to get a big enough support; but nevertheless $L_j q_{n_j}$
is still a sequence of rigidity times.

Now we make the argument more precise.

\vskip 3mm
{\it Support of $\Phi_\beta^{(L_j q_{n_j})}$}

Let $n =n_j$ be in $\cal S$. If we assume for concreteness $n$ odd and $\varepsilon_n(\beta) = +1$, using equations (\ref{disczer}) to (\ref{discbeta2}),
we obtain that the value of $\Phi_\beta^{(q_n)}$ is
\begin{eqnarray*}
-1 \text{ on } &&I_{k,n}^0 := [{k \over q_n} + {u_n(k, 0)  \|q_n \alpha\| \over q_n}, \  {k \over q_n} +  {\|q_n \beta\| \over q_n}
+ {u_n(k, \beta)  \|q_n \alpha\| \over q_n}[, \\
1 \text{ on } &&I_{k,n}^1 := [{k \over q_n} +  {1 \over 2q_n} + {u_n(k, \frac12)  \|q_n \alpha\| \over q_n},
\ {k \over q_n} + {1 \over 2q_n}+  {\|q_n \beta\| \over q_n}  +  {u_n'(k, \frac12 + \beta)  \|q_n \alpha\| \over q_n}[,
\end{eqnarray*}
and 0 elsewhere.

Since $u_n(k, 0)  \|q_n \alpha\|$ and $u_n(k, \beta)  \|q_n \alpha\|$ are both $\leq q_n \|q_n \alpha\|$ and
$\|q_n \beta\| \geq c_n(\beta) q_n  \|q_n \alpha\| \geq  2 q_n  \|q_n \alpha\|$
for $n \in \cal S$, we have $\|q_n \beta\| +u_n(k, \beta)  \|q_n \alpha\| - u_n(0, \beta)  \|q_n \alpha\| > q_n  \|q_n \alpha\|$. Hence
the length of  $I_{k,n}^0$ is bigger than $\|q_n \alpha\| \geq \frac12 q_n^{-1}  a_{n+1}^{-1}$. Similarly, the length of  $I_{k,n}^1$ is bigger than
$\|q_n \alpha\|$.

Let $L_j = [\delta \, a_{n_j+1}]$, where $0 < \delta < 1/2$ be a constant. Therefore  $L_j \|q_{n_j} \alpha\| < \delta /q_{n_j}$
(so that $R_\alpha^{L_j q_{n_j}} \to Id$ and $\Phi_\beta^{(L_j q_{n_j})}$ tends to 0 in measure). Let us consider
the sum $\Phi_\beta^{(L_j q_{n_j})}$.

The measure of the subset $J^0(k,n_j)$ (resp. $J^1(k,n_j)$) of $[{k \over q_{n_j}}, \ {k+1 \over q_{n_j}}[$ on which
$\Phi_\beta^{(L_j q_{n_j})} \leq -1$ (resp. $\Phi_\beta^{(L_j q_{n_j})} \geq 1$)
is the measure of the union $\widetilde I_{k,n_j}^0$ (resp. $\widetilde I_{k,n_j}^1$) of the intervals translated of
$I_{k,n_j}^0$  (resp. of $I_{k,n_j}^1$) by $u \theta_n$, $u = 0, ...,  L_j-1$; therefore is bigger than $\delta a_{n_j+1} \|q_{n_j} \alpha\|$.

Therefore, the measure of the union $A_{n_j}^0 = \bigcup_{k=0}^{q_{n_j}-1} \, J^0(k,n_j)$ is bigger than
$$q_{n_j} \, L_{n_j} \|q_{n_j} \alpha\| = \delta a_{n_j+1} \, q_{n_j} \|q_{n_j} \alpha\| > \delta a_{n_j+1}  \, q_{n_j}/ (a_{n_j+1}  \, q_{n_j} +q_{n_j-1})
\geq \frac12 \delta.$$

Hence, along the sequence $\cal S$, $\Phi_\beta^{(L_j q_{n})}$ does not tend to 0 in measure and equation (\ref{funeq1}) has no measurable solution.
\eop

\begin{remark} \  \label{noteven1} \em Suppose that there are infinitely many even denominators, but that the folowing condition is satisfied:
\begin{eqnarray}
M:=\sup_{n \, : \, q_n \text{ even }} a_{q_n+1} < +\infty. \label{majEvn1}
\end{eqnarray}

Then, if $\lim_n \|q_n \beta \| = 0$, the same proof as in case B) above applies. Indeed, it suffices to check that, for
the subsequence $(n_j)$ such that $c_{n_j}(\beta) > \frac14$ given by Lemma \ref{qnbeta0}, $q_{n_j}$ is odd for $j$ big enough.

Suppose that $q_{n_j}$ is even. Then, by \eqref{f_3b} and \eqref{majEvn1}, we have:
$$\|q_{n_j} \beta \| \geq \frac14 q_{n_j} \|q_{n_j} \alpha\| \geq \frac14 {1 \over a_{n_j +1} +2} \geq \frac14 {1 \over M+2}.$$
Therefore, $q_{n_j}$ is odd, once $j$ satisfies $\|q_{n_j} \beta \| <  \frac14 {1 \over M+2}$.
\end{remark}

\begin{remark} \  \label{nonreg1} \em In the next section the existence of values of $\beta$ giving non-regularity will be shown. Regularity or non-regularity of the cocycle
is related to the behavior of the sequence $(c_{n_{j}}(\beta))_{j \geq 1}$. (Recall that $c_n(\beta)$ was defined in~\eqref{defcn}.)

Suppose that for a subsequence $(n_{j_\ell})$ and a finite constant $K$, we have
$K^{-1} \leq c_{n_{j_\ell}}(\beta) \leq K$, then the skew map
$R_{\alpha, \Phi_\beta}$ is regular. Indeed, the overlapping of the support occurs for at most $K+1$ translation by $q_{n_\ell} \alpha$.
This implies that on sets with a measure bounded away from 0, $\Phi_\beta^{(L_j q_{n_{j_\ell}})}$ takes a fixed non zero integer value, which therefore is an essential value
of the cocycle. So we get that the cocycle is regular.

But $\|q_{n_j} \beta\|/ q_{n_j} = c_{n_j}(\beta) \|q_{n_j} \alpha\|$ can be much bigger than $\|q_{n_j} \alpha\|$, in which case there is a big
overlapping of the translates of $I$. Therefore, if  $ c_{n_j} \uparrow \infty$, non-regularity can occur. In the next section we will see that this can be effectively the case.
\end{remark}

\vskip 3mm
{\bf Proof of Theorem~\ref{Mainthm} part 2c)}

This part now follows directly from Theorem~\ref{thmfunctEqu1}.

\vskip 3mm
{\bf Proof of Theorem~\ref{Mainthm} part 2d)}

This part will follow Theorem \ref{evenThm1} below. We need some preliminary results.

Let $J = \{0, \beta, \frac12, \beta-\frac12\}$. For $\beta \not \in \Z /2$, for every $N \geq 1$, the set of discontinuities of the ergodic sum $\Phi_\beta^{(N)}$ 
consists of the distinct points $t- j\alpha {\rm \mod 1}$, where $0 \leq j < N$ and $t \in J$, with jumps $\pm 1$. We write $\{0= \gamma_{N,1} < ... <\gamma_{N,4N}\}$
for the elements in this set listed in natural order.

By minimality of the rotation $R_\alpha$, the following lemma holds:
\begin{lem} \label{rigidLem} Let $\lambda > 0$. For every $\varepsilon >0$, there is $L(\varepsilon)$ such that, for any $L \geq L(\varepsilon)$, we have
$\|b \,\alpha\| \leq \varepsilon$, for some  $b \in [\lambda \, L, 2 \lambda \, L]$. 
\end{lem}

\begin{proposition} \label{rational1} Let $\alpha$ be an irrational number such that the sequence $(a_n)$ of its partial quotient does not tend to infinity. Then,
if $\beta \in (\Q \alpha + \Q) \setminus (\Z\alpha + \Z)$, $\Phi_\beta$ is not a measurable coboundary.
\end{proposition}
\proof We have $\beta = {\ell \over s} \alpha + {r\over s}$, with $\ell, r, s$ integers and $s \not = 0$.
By assumption, there is a fixed integer $a$ and subsequence $(n_k)$ such that $a_{n_k+1} = a$.

Let $n$ be such that $a_{n+1} =a$. We have (cf. (\ref{f_4})) $\|k \alpha \| \geq \|q_n \alpha \|, \forall k \in [0, q_{n+1}[$, and
$$\|q_n \alpha \| \geq {1\over q_n+q_{n+1}} = {1\over (a_{n+1}+1)q_n+q_{n-1}} \geq {1\over 2+a} {1\over q_n}.$$

Put $R = 2s$. Let $t \in J$. By the previous inequalities, we have for an integer $\ell(t)$, for $k \leq q_n/R$:
\begin{eqnarray}
\|k \alpha - t\| \geq {1\over R} \|(R k + \ell(t))\alpha\| \geq {c \over q_n}, \text{ with } c :={1\over R} {1\over 2 + a}, \label{minor1}
\end{eqnarray}
since $Rk + \ell(t) < q_n+q_{n-1} \leq a_{n+1}q_n+q_{n-1} = q_{n+1}$, for $n$ big enough and $k \leq q_n/R$.

If $(\varepsilon_i)$ is a sequence of positive numbers tending to 0, by Lemma \ref{rigidLem}
we can choose $(k_i)$ and $b_i \in [{q_{n_{k_i}} \over 2R}, {q_{n_{k_i}} \over R}]$ such $\|b_{i}  \alpha\| \leq \varepsilon_i$.
So $(b_i)$ is a sequence of rigidity times for the rotation by $\alpha$.
This implies that, if $\Phi_\beta$ is a measurable coboundary, the ergodic sums $\Phi_\beta^{(b_i)}$ tends to 0 in measure.

On the other hand, for $\ell = 1, ..., 4b_i-1$, $\Phi_\beta^{(b_i)}$ is constant on the intervals $]\gamma_{b_i,\ell}, \gamma_{b_i,\ell+1}[$ and we have
$\gamma_{b_i, \ell+1} - \gamma_{b_i,\ell}   \geq {c \over q_{n_{k_i}}} \geq {c \over 2 R b_i}$ by~(\ref{minor1}).
Therefore $|\Phi_\beta^{(b_i)}| \geq 1$ on a set of measure bounded away from 0. This gives a contradiction.
\eop

\begin{lem} \label{minorigidlem} Let $n \geq 1$ and $\delta > 0$ be such that $\|q_n \beta\| \geq \delta$. 
For $b \in [{\delta \over 4} \, q_n, \ {\delta \over 2} \, q_n]$, we have
\begin{eqnarray}
\|\beta - j\alpha\| \geq  {\delta^2 \over 8}  \, {1 \over b}, \text{ for } |j| < b. \label{minorigidIneq}
\end{eqnarray}
\end{lem}
\proof From 
$\|q_n(\beta - j \alpha)\| \geq \|q_n \, \beta \| - \|q_n \, j \alpha\| \geq \delta - |j| \|q_n \, \alpha\| \geq \delta - b / q_n \geq {\delta \over 2}$, it follows:
$$\|\beta - j\alpha\| \geq q_n^{-1} \, \|q_n(\beta - j \alpha)\| \geq {\delta \over 2} \, q_n^{-1} \geq {\delta \over 2}  \, {\delta \over 4}  \, b^{-1}
= {\delta^2 \over 8} \, b^{-1}, \text{ for } 0 \leq j < b. \eop$$

\begin{proposition} \label{regularprop} Let $\beta \in ]0, \, 1[$. If there is a subsequence $(q_{n_k})_{k \geq 1}$ of odd denominators of $\alpha$
such that, for some $\delta \in ]0, \, \frac12[$, $\|q_{n_k} \beta\| \geq \delta, \, \forall k \geq 1$, then $\Phi_\beta$ is not a measurable coboundary 
for the rotation $R_\alpha$.
\end{proposition}
\proof Let $(\varepsilon_i)$ be a sequence of positive numbers tending to 0. By Lemma \ref{rigidLem} we can choose $(k_i)$ 
and $b_i \in [{\delta \over 4} \, q_{n_{k_i}}, \ {\delta \over 2} \, q_{n_{k_i}}]$ such $\|b_{i}  \alpha\| \leq \varepsilon_i$. 
Suppose that $\Phi_\beta$ is a measurable coboundary. Then the ergodic sums $\Phi_\beta^{(b_i)}$ tends to 0 in measure.
We will show that this is not possible.  

Since $q_{n_k}$ is odd, we have $\|q_{n_k}  \frac12\| = \frac12 \geq \delta$. 
By Lemma \ref{minorigidlem}, for the ergodic sum $\Phi_\beta^{(b_i)}$  the discontinuities of type 0 are ``well separated'' from the discontinuities of type $\beta$ 
and of type $\frac12$, since we have, with $c = {\delta^2 \over 8}$,
$$\inf_{|\ell| < b_i} \, (\|\beta - \ell \alpha\|, \, \|\frac12 - \ell \alpha\|) \geq {c \over b_i}.$$

Denote by $\gamma_0$ any discontinuity of type 0. Let $\gamma^-$ (resp. $\gamma^+$) be the nearest discontinuity of type $\beta$ or $\frac12$ 
at left (resp. at right) of $\gamma_0$. 

The possible jumps between $\gamma^-$ and $\gamma_0$ (resp. $\gamma_0$ and $\gamma^+$) are only $-1$. Therefore,  $\Phi_\beta^{(b_i)}$ is non increasing 
on  $]\gamma^-, \  \gamma^+[$ and moreover its value is decreased by $-1$ at the point $\gamma_0$. It follows that
$$\Phi_\beta^{(b_i)}(x) \geq 1, \text{ for } x \in ]\gamma^-, \  \gamma_0[, \ \text{ or } \ \Phi_\beta^{(b_i)}(x) \leq -1, \text{ for } x \in ]\gamma_0, \  \gamma^+[.$$

As we know that the distance between $\gamma^-$ and $\gamma_0$ (resp. $\gamma_0$ and $\gamma^+$) is $\geq c/ b_i$, we conclude that, on the whole circle,
$|\Phi_\beta^{(b_i)}| \geq 1$ on a set of measure bounded away from 0. We get a contradiction.\eop

\vskip 3mm
\begin{thm}\label{evenThm1} Let $\alpha$ be an irrational number such that there are finitely many even denominators or 
$\sup_{n \, : \, q_n \text{ even }} a_{n+1} < +\infty$. Then, for $\beta \not \in \Z\alpha + \Z$, equation
\begin{eqnarray}
\Phi_\beta(x) = g(x+\alpha) - g(x), \text{ for }\mu-\text{a.e. }x\in\T,  \label{funeq2}
\end{eqnarray}
has no measurable solution~$g$.
\end{thm}
\proof If $\limsup_{n \, : \, q_n \text{ odd }} \|q_{n} \beta\| > 0$, the result follows  from Proposition~\ref{regularprop}. Therefore we can assume
$\lim_{n \, : \, q_n \text{ odd }} \|q_{n} \beta\| = 0$.

Observe that, if there is $n_0$ such that all denominators $q_n$ are odd for $n \geq n_0$, then the result follows from Theorem~\ref{thmfunctEqu1}.
Let us consider now the case where there are infinitely many even denominators and let denote by $q_{n_1}<q_{n_2}<\ldots$ their sequence.

Since the denominators $q_{n_k+1}$ and $q_{n_k-1}$ are odd, it holds $\lim \|q_{n_k+1} \beta \| = \lim \|q_{n_k-1} \beta \| = 0$.
From the relation $q_{n_k+1} \beta = a_{n_k+1} q_{n_k} \beta + q_{n_k-1}\beta$, it follows $\lim_k  \|a_{n_k+1} q_{n_k} \beta\|= 0$.

Let $R$ denote the integer $R=\prod_{j \in J} A_j$, where $\{A_j:\: j \in J\}$ is the finite set of values taken by the $a_{n_k + 1}$'s.
We get $\lim_k \|q_{n_k} R \beta\| = 0$. Therefore, for the whole sequence $(q_n$), we have $\lim_n \|q_{n} R \beta\| = 0$.

If $R\beta \not \in \Z \alpha + \Z$, then part B) of the proof of Theorem~\ref{thmfunctEqu1} applies to $R \beta$
(see Remark \ref{noteven1}). This shows that the function $\Phi_{R\beta}$, and so as well $\Phi_{\beta}$, is not a measurable coboundary.

Finally the remaining case $\beta \in (\Q \alpha + \Q) \setminus (\Z\alpha + \Z)$  is treated in Proposition \ref{rational1}. \eop

Remark that, if the hypothesis of the theorem is not satisfied, then the situation is that of Theorem \ref{Mainthm} 1b) and equation (\ref{funeq2}) has a solution 
for uncountably many $\beta$'s.

\section{\bf Ostrowski expansion and non-regularity for an exceptional set} \label{Except}

\

As remarked above, the proof in case B) of  Theorem~\ref{thmfunctEqu1} gives a result weaker than ergodicity. Actually, we will show that
in that case there is a set of values of $\beta$ for which ergodicity (as $\Z$-valued cocycle) fails and the cocycle $\Phi_\beta$ is non-regular.

Denote by ${\cal U}(\T)$ the group of measurable functions from $\T=[0,1[$ to the group ${\cal U}$ of complex numbers of modulus~1.

The non-regularity result is based on the following observation (cf.~\cite{Co09}): if $g$ is cohomologous to $g_1$ and to $g_2$, two functions with values respectively
in closed subgroups with an intersection reduced to $\{0\}$,
then ${\cal E}(g) = \{0\}$. This implies:
\begin{lem}\label{sgroupdis} If $\varphi$ is a $\Z$-valued cocycle such that there exists $s \not \in \Q$ for which the multiplicative
equation $e^{2\pi i s \varphi} =  \psi \circ  R_{\alpha}/\psi$ has a measurable solution $\psi:\T\to{\cal U}$, then ${\cal E}(\varphi) = \{0\}$. If
$\varphi$ is not a coboundary, then  ${\overline {\cal E}}(\varphi) = \{0, \infty \}$ and $\varphi$ is non-regular.
\end{lem}

Let us consider the function $\psi_{\beta,s}:= e^{2\pi i s 1_{[0,\beta]}}$ on the circle and the multiplicative functional equation
\begin{equation}\label{equafonc1}
e^{2\pi is 1_{[0,\beta]}} = e^{2\pi it} \, {R_\alpha f / f},
\text{ where } (\beta, s, t)\in [0,1[\times \R \times \R \text{ and } f \in {\cal U}(X),
\end{equation}

This equation was studied by W.\ Veech in \cite{Ve69}, then by K.\ Merril \cite{Me85} who gave a sufficient condition on $(\beta,
s, t)$ for the existence of a solution, then by M.\ Gu\'enais and F.\ Parreau \cite{GuPa06} who gave a necessary and sufficient condition for~(\ref{equafonc1}) to have a measurable solution
and extended it to more general step functions. The conditions are expressed in terms of the so-called Ostrowski expansion of a real $\beta$.
For $r \geq 1$, we put
\begin{eqnarray*}
&&H_r(\alpha):=\left\{\sum_{n\geq 0}b_n q_n\alpha \mod 1 ,\  (b_n)_n \in \M Z^{\M N},
\text{ such that } \sum_{n\geq 0} \left({|b_n| \over a_{n+1}}\right)^r <+\infty \right\},\\
&&H_\infty(\alpha):= \left\{\sum_{n\geq 0} b_n q_n\alpha \mod 1,\ {|b_n| \over a_{n+1}} \rightarrow 0\right\}.
\end{eqnarray*}

Recall the following characterization (\cite{GuPa06}):
\begin{eqnarray*}
&&H_r(\alpha)=\{\beta\in\R: \, \sum_{n\geq0} \|q_n \beta\|^r <+\infty\},\ \ H_\infty(\alpha)=\{\beta\in\R: \, \|q_n \beta \| \to 0\}.
\end{eqnarray*}
When $\alpha$ is not of bounded type, $H_\infty(\alpha)$ is an uncountable additive subgroup of $\R$.

\begin{thm} (\cite{GuPa06}) \label{NSCdtion2} Equation (\ref{equafonc1}) has a solution $f \in {\cal U}(X)$ for the parameters $(\beta, s, t)$
if and only if there is a sequence $(b_{n})$ in $\Z$ such that:
\begin{eqnarray}
&&\beta=\sum_{n\geq 0}b_{n}q_{n}\alpha \textrm{ mod }1, \text{ with } \sum_{n\geq 0}\frac{|b_{n}|}{a_{n+1}} = C_1 <\infty, \label{OstroCondi}\\
&&\sum_{n\geq 0}\Vert b_{n}s\Vert^{2}<\infty, \ \ t=k \alpha-\sum_{n\geq 0} [b_{n}s] \, q_{n}\alpha \textrm{ mod }1, \text{ for an integer } k.
\end{eqnarray}
\end{thm}

The size of $c_n(\beta)$, a key point in the proof of Theorem \ref{thmfunctEqu1}, is related to the $b_n$'s in the expansion of $\beta$.
For a non trivial $\beta$, when the $b_n$'s are bounded, it can be shown by the method of Theorem \ref{thmfunctEqu1} that $\Phi_\beta$ is ergodic.
At the opposite, a fast growth of the sequence $(b_n)_{n \geq 1}$ implies the non-regularity of the cocycle:

\begin{thm}\label{val-ess}  Let $R_\alpha$ be the rotation by $\alpha$ with unbounded partial quotients. If $\beta \not \in \Z \alpha + \Z$
satisfies (\ref{OstroCondi}) with the lacunarity condition  $\sum_n (b_n /b_{n+1})^2 < \infty$, then $\Phi_\beta = \varphi_{\beta, \frac12}$ defines
a non-regular cocycle (and therefore the skew product $R_{\alpha, \Phi_\beta}$ is not ergodic).
\end{thm}
\Proof By Theorem \ref{NSCdtion2}, if $\beta$ satisfies (\ref{OstroCondi}), for $s$ in the set $\{s: \sum_{n\geq 0}\Vert b_{n}s\Vert^{2}<\infty\}$,
there is a solution of (\ref{equafonc1}). Moreover, the set of such $s$ is uncountable if
$\sum_n (b_n /b_{n+1})^2 < \infty$. There are thus $\beta \not \in \alpha \Z + \Z$, $s \not \in \Q$, $t \in \R$
and $\psi \in {\cal U}(X)$ of modulus 1 such that $e^{2\pi i s 1_{[0,\beta]}} = e^{2\pi it} \, \psi \circ R_\alpha / \psi$.

For this choice of $(\beta, s)$, $e^{2\pi i s(1_{[0,\beta]} - 1_{[0,\beta]}\circ R_\frac12)}$ is a multiplicative coboundary.
On the other hand, we have shown that $1_{[0,\beta]} - 1_{[0,\beta]}\circ R_\frac12 = \varphi_{\beta, \frac12} = \Phi_\beta$
is not an additive coboundary. Lemma \ref{sgroupdis} shows that ${\overline {\mathcal E}}(\varphi_{\beta, \frac12})= \{0, \infty \}$,
which implies the non-regularity of $\Phi_\beta$.\eop

{\it Remark}: The previous result gives an explicit value of $\gamma$, namely $\gamma=\frac12$, such that $\varphi_{\beta, \gamma}$ is non regular.
A generic result also holds  (cf.\ \cite{Co09}): if $\beta$ satisfies (\ref{OstroCondi}) with the lacunarity condition
$\sum_n (b_n /b_{n+1})^2 < \infty$, then, for a.e. $\gamma$, $\varphi_{\beta, \gamma}$ is a non-regular cocycle.

The previous result is for $\alpha$ of Liouville type. At the opposite, if we take $\alpha$ with bounded partial quotients,
as we have seen (cf. Proposition~\ref{zergodic}), for $\beta \not \in \Z \alpha + \Z$, the $\Z$-valued cocycle $\Phi_\beta$ is ergodic.

\begin{rem} Let us consider $\varphi_{\beta,\gamma}$ (cf. Notation \ref{phibetaj}). For $\beta, \gamma \in \, ]0,1[$, with $\beta + \gamma < 1$,
this step function reads $1_{[0, \beta[} - 1_{[\gamma, \beta + \gamma[}$. Its Fourier coefficients of $\varphi_{\beta,\gamma}$ are
${1 \over 2 \pi i n} (e^{-2\pi in \beta} - 1) (-e^{2\pi i n \gamma} - 1)$.

The condition for $\varphi_{\beta,\gamma}$ to be a coboundary with a transfer function
in $L^2(\T)$, i.e., such that the functional equation $\varphi_{\beta,\gamma}= R_\alpha g - g$ has a solution $g$ in$L^2$, is
\begin{eqnarray}
\sum_{n \not = 0} {1 \over n^2} { \|n \beta\|^2 \|n \gamma\|^2 \over \|n \alpha\|^2} < \infty. \label{betalambda}
\end{eqnarray}
The following sufficient condition for the existence of an $L^2$-solution of the coboundary equation have been given is proved in \cite{CoMa14}:
If $\beta, \gamma$ are in $H_4(\alpha)$, then (\ref{betalambda}) holds and there is $g$ in $L^2(\T)$ solution of $\varphi_{\beta,\gamma} =  R_\alpha g - g$.

Therefore, if $\alpha$ has unbounded partial quotients, there is an uncountable set of pairs of real
numbers $\beta$ and $\gamma$ such that $\varphi_{\beta,\gamma}$ is a coboundary $R_\alpha g - g$ for $R_\alpha$ with $g$ in $L^2$.
\end{rem}

\vskip 3mm
\section{\bf Questions}

{\bf Question 1.} Is there a special measure-theoretic property that permits to single out the elements $W=\widetilde{S_g}$ from the $C^{\rm lift}(T^f)$?
For example, is it true that if $S\circ T^k$ has entropy zero for each $k\in\Z$, then so is the entropy of $W$?

{\bf Question 2.} (cf.\ Remark~\ref{uw1}) Given a flow $(R_t)$ on $\zdr$, for each measurable subgroup $G\subset C((R_t)_{t\in\R})$, can we find a special representation
$T^f$ of $(R_t)$ such that $C^{\rm lift}(T^f)$ ``realizes'' $G$? (I.e.,\ a measure-theoretic isomorphism $I$ between the flow and its special representation yields
$I(G)=C^{\rm lift}(T^f)$.)

In particular, does there exist a flow $(R_t)_{t\in\R}$ such that for no special representation $T^f$ of it we have $C(T^f)=C^{\rm lift}(T^f)$?

{\bf Question 3.} (cf.\ Remark~\ref{uw1} and Question~2) Can we find  $\alpha$ and $f$ regular for which $C^{\rm lift}(T^f)$ is not closed?

{\bf Question 4.} Assume that $Tx=x+\alpha$ and $f:\T\to\R^+$ is smooth (we recall that then $T^f$ is rigid). Is it true that $C^{\rm lift}(T^f)=C(T^f)$?

\section{\bf Appendix. Centralizer for uniformly rigid special flows}
\subsection{\bf Continuous centralizer of uniformly rigid flows}

\

Let $(X,d)$ be a compact metric space and let $\ct=(T_t)_{t\in\R}$ be a continuous flow on it, i.e.,\ it is a one-parameter group of homeomorphisms of $X$:
$T_t\in{\rm Homeo}(X)$ for $t\in\R$ and
\begin{equation}\label{app0}
\mbox{the map $(x,t)\mapsto T_tx$ is continuous.\footnote{We also assume that the map $t\mapsto T_tx$ is 1-1 for each $x\in X$.}}
\end{equation}
We then have
\begin{equation}\label{app1}
\mbox{the map $t\mapsto T_t$ is continuous,}\end{equation}
where on ${\rm Homeo}(X)$ we consider the uniform topology: $\rho(V,W):=\sup_{x\in X}(d(Vx,Wx)+d(V^{-1}x,W^{-1}x)$ whenever $V,W\in {\rm Homeo}(X)$
(with this topology ${\rm Homeo}(X)$ becomes a Polish group). Indeed, we only need to show that, whenever $\epsilon>0$, we have $d(x,T_tx)<\epsilon$
for all $x\in X$ and $|t|<\delta$ for some $\delta>0$ which results immediately from the uniform continuity of the map $(x,t)\mapsto T_tx$
on $X\times [-1,1]$.

A flow $\ct$ is called {\em uniformly rigid} if for some sequence $s_n\to\infty$, we have $T_{s_n}\to Id$ uniformly. We can now repeat the
``measurable'' proof from~\cite{Ka-Le} in the continuous setting.

\begin{Prop}\label{pot3} Assume that a flow $\ct=(T_t)_{t\in\R}$ is uniformly rigid. Then the essential topological centralizer
$C^{\rm top}(\ct)/\{T_t:\:t\in\R\}$  is uncountable.
\end{Prop}
\proof Consider
$$H:=\{T_t:\:t\in\R\}\subset \overline{\{T_t:\:t\in\R\}}:=G\subset{\rm Homeo}(X),$$
where $G$ is a Polish group. If $H$ is a proper subgroup, then it must be a set of first category, and hence, it cannot have only countably many cosets
(as $G$  is Polish without isolated points). If $H=G$, then $H$ itself is Polish, and by~\eqref{app1} the map  $t\mapsto T_t$ is continuous.
Since this map is 1-1, by the open map theorem for topological groups, the map $t\mapsto T_t$ has to be a homeomorphism, and the continuity
of the inverse yields a contradiction with the uniform rigidity of the special flow $\ct$.
\eop

\vskip 3mm
{\bf Continuous special flows}

Let $(X,d_X)$ be a compact metric space and $f:X\to\R^+$ continuous. In particular, for some $\eta>0$ we have $f(x)\geq \eta$ for each $x\in X$.
Set $\overline{X}^f=\{(x,r)\in X\times\R:\:0\leq r\leq f(x)\}$. Then $\overline{X}^f$ is a compact metric space with the product metric $d$ (the product of $d_X$
and the Euclidean metric $d_{\R}$  on  $\R$). Let $T:X\to X$ be a homeomorphism.
Define the equivalence relation $\sim$ on $\overline{X}^f$ with the only non-trivial gluing $(x,f(x))\sim (Tx,0)$.  The resulting space denoted by $X^f$ is Hausdorff
and compact (and we could  identify it with $\{(x,r):\:x\in X, 0\leq r<f(x)\}$).
Let $D$ be the quotient metric defined by $D((x,r),(x',r')):=$
\beq\label{metilo}
\\ \Big \{\begin{array}{l}\inf\{d((x,r),(x_1,r_1))+d((x'_1,r'_1),(x_2,r_2))+...+\ d((x'_n,r'_n),(x',r')),\\
\text{ where }(x_i,r_i)\sim(x'_{i},r'_{i}),\ i=1,\ldots,n\}.
\end{array}\eeq
Then  $T^f$ becomes a continuous flow on the compact metric space $X^f$.

\vskip 3mm
{\bf Uniform rigidity of special flows}

\begin{Prop}\label{pot2}  Let $T$ be uniformly rigid, that is, for some increasing sequence $(q_n)\subset\N$ we have $T^{q_n}\to Id$ uniformly.
If there exists $(s_n)\subset\R$  such that $f^{(q_n)}(\cdot)-s_n\to 0$ uniformly, then $T^f_{s_n}\to Id$ uniformly.
\end{Prop}
\proof
For each $(x,r)\in X^f$, we have
$$
D(T^f_{s_n}(x,r),(x,r))=D(T^{f}_{s_n-f^{(q_n)}(x)}T^f_rT^f_{f^{(q_n)}(x)}(x,0),(x,r))=
$$
$$
D(T^f_{s_n-f^{(q_n)}(x)}T^f_r(T^{q_n}x,0),(x,r))=D(T^f_{s_n-f^{(q_n)}(x)}
(T^{q_n}x,r),(x,r))\leq$$
$$
D(T^f_{s_n-f^{(q_n)}(x)}(T^{q_n}x,r),(T^{q_n}x,r))+D((T^{q_n}x,r),(x,r))
$$
and the two last summands are small by~\eqref{app1} (if $n$ is sufficiently large) and the definition of $D$.\eop

Directly from Proposition~\ref{pot3}, we obtain the following.

\begin{cor}\label{capp1}
Under the assumptions of Proposition~\ref{pot2}, the essential (topological) centralizer $C^{\rm top}(X^f,T^f)/\{T^f_t:\:t\in\R\}$ is uncountable.
\end{cor}

\subsection{\bf Smooth special flows over irrational rotations}

\

Let us come back to special flows over irrational rotations ($X=\T$, $Tx=x+\alpha$). Let $f:\T\to\R^+$. For simplicity, we assume that $\int_{\T} f\,d\mu=1$ and set $f_0:=f-1$.
Then, it follows from \cite{He} that if $f$ is absolutely continuous (AC), then $f^{(q_n)}(\cdot)-q_n\to 0$ uniformly. Hence $T^f$ is uniformly rigid.

\begin{cor}\label{pot4}  Let $Tx=x+\alpha$ and $f:X\to\R^+$ be AC. Then the essential topological centralizer of
$T^f$ is uncountable. Moreover, there exists an uncountable set of $\beta\in\T$ such that the functional equation
\begin{equation}\label{appe1}
f(x+\beta)-f(x)=g(x+\alpha)-g(x)\end{equation}
has a solution in continuous functions $g:X\to\R$.
\end{cor}
\proof  The first part follows from the uniform rigidity and Proposition~\ref{pot3}, the second one is a consequence of the first one
and of the result from \cite{Markley} on the form of homeomorphisms commuting with $T^f$.
\eop

We will now show a different (direct) proof (cf. \cite{Le-Ma}) of the fact that whenever $f$ is AC then we can solve~\eqref{appe1}
for uncountably many $\beta$.

For this aim select a subsequence $(q_{n_k})_{k\geq1}$ of denominators of $\alpha$ so that
\begin{equation}\label{appe2}
\sum_{k\geq1}\|f_0^{(q_{n_k})}\|_{C(\T)}<+\infty\text{ and }\sum_{k\geq1}\|q_{n_k}\alpha\|<+\infty\end{equation}
(remembering that $f_0^{(q_n)}\to0$ uniformly and $\|q_n\alpha\|\to0$).
We have, for each $x\in\T$ and $k\geq1$,
$$f_0^{(q_{n_k})}(x+\alpha)-f_0^{(q_{n_k})}(x)=f_0(x+q_{n_k}\alpha)-f_0(x).$$
By replacing  $x$ by   $x+\sum_{j<k}q_{n_j}\alpha$, we obtain
\begin{eqnarray*}
f_0^{(q_{n_k})}(x+\sum_{j=0}^{k-1}q_{n_j}\alpha +\alpha)- f_0^{(q_{n_k})}(x+\sum_{j=0}^{k-1}q_{n_j}\alpha)=
f_0(x+\sum_{j=0}^{k-1}q_{n_j}\alpha +q_{n_k}\alpha)- f_0(x+\sum_{j=0}^{k-1}q_{n_j}\alpha).
\end{eqnarray*}

Now, the RHS  of the above equality is telescopic, and when we sum it up, by~\eqref{appe2}, we obtain $f_0(x+\beta)-f_0(x)$ with  $\sum_{k\geq 1}q_{n_k}\alpha=\beta$,
while for the LHS the series
$$\sum_{k\geq1} f_0^{(q_{n_k})}(x+\beta_k),\text{ where }   \beta_k=\sum_{j=0}^{k-1}q_{n_j}\alpha,$$
converges uniformly as it converges absolutely by~\eqref{appe2}.
By~\eqref{appe2}, we have
$\sum_{k\geq1} f_0^{(q_{n_k})}(\cdot+\beta_k)=g$.
Hence we obtain~\eqref{appe1}. Note finally that if in the above reasoning we replace
$q_{n_k}$ by $\epsilon_kq_{n_k}$, with $\epsilon\in\{0,1\}^{\N}$ (with infinitely many $k$ for which $\epsilon_k=1$), using a unicity argument in the
Ostrowski expansion of $\beta$, we obtain an uncountable set of $\beta\in\T$ for which we can solve~\eqref{appe1}.

The above method can be also applied  when the roof function
$f=\sum_{n\in\Z}a_ne^{2\pi nx}$ satisfies $a_n={\rm o}(1/|n|)$. Indeed, as proved in \cite{Le-Ma}, under this assumption, $f^{(q_n)}_0\to 0$ in $L^2(\T)$. It follows that the corresponding special flow is rigid, whence its essential centralizer is uncountable. But by repeating the above proof, we obtain:

\begin{proposition}\label{omale}
Let $f=\sum_{n\in\Z}c_ne^{2\pi inx}$ satisfy $c_n={\rm o}(1/|n|)$. Then for each irrational $\alpha$ the set of $\beta$ for which we can solve~\eqref{funeq1}
with $g\in L^2(\T)$ is uncountable. Equivalently, the essential liftable centralizer is uncountable.
\end{proposition}

\begin{proposition}\label{series1}
Let $f$ be in $L^2(\T)$ such that, for an irrational $\alpha$ and a strictly increasing sequence $(r_n)_{n \geq 1}$,
\begin{eqnarray*}
\ \|r_n \alpha\| \to 0, \ \|f^{(r_n)}\|_2 \to 0.
\end{eqnarray*}
Then the set of $\beta$ for which we can solve
$f(x+\beta) - f(x) = g(x+\alpha) - g(x)$ with $g\in L^2(\T)$ is uncountable.
\end{proposition}

\begin{remark}
\em For the smooth case $C^2$, A.\ Kanigowski gave a Fourier analysis type argument showing that the set of $\beta$ for which~\eqref{appe1} can be solved is residual.
\end{remark}

\subsection{\bf Special flow with H\"olderian roof function and trivial liftable centralizer}\label{hoelder}

\

The aim of this section is to show the following result (to be compared with Corollary~\ref{pot4}).

\begin{Prop}\label{p:hoelder} For each $\alpha$ with bounded partial quotients, there is $F$ which is H\"older continuous with any H\"older exponent $0<\kappa<1$ and such that the functional equation
\begin{eqnarray}
F(x+\beta) - F(x) = g(x+\alpha) - g(x) \label{funeqF}
\end{eqnarray}
has a measurable solution $g$ only for $\beta \in \Z \alpha + \Z$. In other words, the liftable centralizer of the special flow $R_\alpha^F$ is trivial.\end{Prop}

To prove Proposition~\ref{p:hoelder}, given $\alpha$ with bounded partial quotients, we will construct below a class of ergodic continuous cocycles $F$
such that the functional equation~(\ref{funeqF}) has a measurable solution $g$ only for $\beta \in \Z \alpha + \Z$. Our construction is similar to the constructions
of ergodic cocycles  using lacunary Fourier series, see Voln\'y \cite{Vo03}, Br\'emont \cite{Br10}.
We start with two remarks.

1) Recall that a sequence $\Lambda = (n_k)$ of positive integers  is called  lacunary if $\inf_k {n_{k+1} \over n_k} > 1$.
We say that $f \in L^1(\T)$ is a {\it lacunary} if $f(x) = \sum_{n \in \Lambda} \, c_n(f) \,  e^{2\pi i n x}$,
where $\Lambda$ is a lacunary sequence.

Recall that if $f$ is lacunary, then, as $f(x+\beta) - f(x) $ is also lacunary, by a result of M.\ Herman (edited in \cite{He04}), the cocycle $f(x+\beta) - f(x)$
is a measurable coboundary if and only if it is a coboundary in $L^2$.

Therefore, if $F$ is lacunary, a measurable solution $g$ of (\ref{funeqF}) exists if and only if
\begin{eqnarray} \label{akjp}
\sum_{n \not = 0} |c_n(F)|^2 {|\sin (\pi n \beta)|^2 \over |\sin (\pi n \alpha)|^2} < \infty.
\end{eqnarray}

2) Let $\alpha$ be an irrational with bounded partial quotients. Then, the sequence $(q_n)$ of denominators of $\alpha$ is lacunary. Indeed, setting $A:= \max_n a_n$,
for all $n \geq 3$, we have: $q_{n-1} \leq A q_{n-2} + q_{n-3} \leq (A + 1) q_{n-2}$;
whence
\begin{eqnarray}
q_n \geq q_{n-1} + q_{n-2} \geq (1 + {1 \over A + 1}) \, q_{n-1}. \label{lacunBpq}
\end{eqnarray}

Moreover, see~\eqref{f_3}, we have: ${q_k \over q_{k+1} + q_k} \leq q_k \|q_k \alpha\| \leq {q_k \over a_{k+1}q_{k}+q_{k-1}} \leq {1 \over a_{k+1}}$, so
\begin{eqnarray}
{1\over A + 1} \leq q_k \|q_k \alpha\| \leq 1. \label{minbpq1}
\end{eqnarray}

\begin{lem} \label{almLacun} For each irrational $\alpha$ and $n\geq1$, we have:
\begin{eqnarray}
q_1 +q_2+ ... + q_n &&\leq 2 q_{n+1}, \label{maj1qn}\\
{1\over q_{n+1}} + ... +  {1\over q_{n+k}} + ... &&\leq \frac{C}{ q_{n+1}}, \label{maj2qn}
\end{eqnarray}
where $C = 5+2\sqrt 5$.
\end{lem}
\proof \ 1) Inequality (\ref{maj1qn}) is clearly satisfied for $n= 0, 1$. If we assume that the inequality is true for $n-1$ and $n$, then:
$q_1 +q_2+ ... + q_{n-1} +q_n + q_{n+1} \leq 2 q_{n} + q_n + q_{n+1} \leq 2 (q_n + q_{n+1}) \leq 2 q_{n+2}$, so~\eqref{maj1qn} holds.

2) For $n \geq 1$ fixed, set $r_0 = q_n$, $r_1 = q_{n+1}$, $r_{k+1} = r_k + r_{k-1}$, for $k \geq 1$.
It follows immediately by induction that $q_{n+k} \geq r_k, \, \forall k \geq 0$.

Denote $c={1\over \sqrt 5}$ and let $\lambda_1= \frac{\sqrt 5}2 + \frac12$, $\lambda_2 = -\frac{\sqrt 5}2 + \frac12$ be the two roots of the polynomial
$\lambda^2 - \lambda -1$. Since $\lambda_j^{\ell+1}=\lambda_j^\ell+\lambda_j^{\ell-1}$ for each $j=1,2$ and $\ell\geq1$, we obtain by induction that:
\begin{eqnarray}
q_{n+k} \geq r_k = c \, \lambda_1^k \, (q_{n+1} - \lambda_2 q_n) - c \, \lambda_2^k \, (q_{n+1} - \lambda_1 q_n),
\ k \geq 0, n \geq 1. \label{Fibo1}
\end{eqnarray}
Take $k\geq1$.
Since $\lambda_2<0$ and $\left|1-\frac{\lambda_1q_n}{q_{n+1}}\right|<\lambda_1$ (as $\lambda_1>1$), from~(\ref{Fibo1}), we obtain
$$q_{n+k} \geq c \, \lambda_1^k \,\left(1 - \lambda_1 \, \left({|\lambda_2|\over \lambda_1}\right)^k\right) \, q_{n+1}=
c\lambda_1^kq_{n+1}\left(1-|\lambda_2|\left(\frac{|\lambda_2|} {\lambda_1}\right)^{k-1}\right)\geq
c\lambda_1^kq_{n+1}(1-|\lambda_2|).$$
It follows that for $k\geq1$, we have
$q_{n+k} \geq c_1 \, \lambda_1^k \, q_{n+1}$, with $c_1 :=\frac{3\sqrt5-5}{10}$. Finally, we obtain
\begin{eqnarray*}
 q_{n+1} \, \sum_{k \geq 1} q_{n+k}^{-1}
 \leq c_1^{-1} \, \sum_{k \geq 1}  \lambda_1^{-k} =c_1^{-1} \, (\lambda_1 - 1)^{-1}
 = 5+2\sqrt 5. \eop
\end{eqnarray*}

Let $s=(m_k)$ be an increasing sequence of positive integers and $\delta >0$. We set
\begin{eqnarray}
F_1(x) =  \sum_{k \geq 1} \, {\sin(2\pi q_{k}x) \over q_{k}}, \ F_s(x) =  \sum_{k \geq 1} \, {\sin(2\pi q_{m_k}x) \over q_{m_k}},  \ F = F_s + \delta F_1.
\label{defF1}
\end{eqnarray}
\begin{prop} Let $\alpha$ be such that the sequence $(q_n)$ is lacunary (in particular, we can take $\alpha$ with bounded partial quotients). If $\beta$ is such that equation~(\ref{funeqF}) for $F$ has a measurable
solution, then $\beta \in \Z \alpha + \Z$.
The function $F$ satisfies the regularity condition:
\begin{eqnarray}
|F(x+h) - F(x)| \leq C|h| \, \log ({1 \over |h|}). \label{regFm}
\end{eqnarray}
In particular, $F$ is H\"olderian with any exponent $0<\kappa<1$.

Moreover, if $\alpha$ has bounded partial quotients, then the sequence $s$ and $\delta$ can be chosen so that  the extension map $R_{\alpha, F}$ on $\T \times \R : (x,y) \to (x+\alpha, y + F(x))$
is ergodic.
\end{prop}
\proof \ 1)  Since, by assumption, the sequence $(q_k)$ of denominators of $\alpha$ is lacunary,  the function $F_s+\delta F_1$ is lacunary and so is the function
$(R_\beta - I) \, (F_s + \delta F_1)$. It follows by \cite{He04} that
if $(R_\beta - I) \, (F_s + \delta F_1)$ is a measurable coboundary, then  equation~(\ref{funeqF}) can be solved in $L^2$, which (by~\eqref{akjp}) implies:
$$\delta^2 \, \sum_{j \not \in s} {1 \over q_{j}^2} {\|q_j\beta\|^2  \over \|q_j \alpha\|^2} +(1 + \delta)^2
\sum_{k}  \, {1 \over q_{m_k}^2} {\|q_{m_k} \beta\|^2  \over \|q_{m_k} \alpha\|^2} < \infty.$$
It follows that $\sum_{k} {1 \over q_{k}^2} {\|q_k\beta\|^2  \over \|q_k \alpha\|^2} < \infty$, which implies that there is $k_0$ such that
$\|q_k\beta\| \leq \frac14  \, q_{k} \, \|q_k \alpha\|$, for $ k \geq k_0$ and therefore, by Lemma \ref{qnbeta0}, $\beta \in \Z\alpha + \Z$.

2) By Lemma \ref{almLacun}, for any $L\geq1$, we have: $|F(x+h) - F(x)|$
\begin{eqnarray*}
&&\leq C' \sum_{k=1}^{L-1} \, {|\sin(2\pi q_{k} (x+h)) - \sin(2\pi q_{k}x)| \over q_{k}} + 2 \sum_{k \geq L} \, {1 \over q_{k}} \leq C' |h|L
+ {C' \over q_{L}}
\end{eqnarray*}
for a constant $C'>0$.
Recall that $q_n \geq C' \lambda_1^n$, with $\lambda_1 > 1$ (cf.\ the proof of Lemma \ref{almLacun}).  It suffices to show (\ref{regFm}) for $h \in ]0, 1/\lambda_1]$. Since $h\leq \lambda_1^{-1}$, there exists $y=y(h) \geq 1$ such that $h y = \lambda_1^{-y}$. We have $ \lambda_1^{-y} \geq h$, whence $y \leq {1 \over \ln \lambda_1} \, \ln {1 \over h}$.

For $0\leq h \leq e^{-\frac12}$, we have $1 \leq 2 \ln {1 \over h} $. Let us take $L = [y]+1$ (that is, $L=L(h)$). We have
$$
hL + {1 \over q_{L}}={\rm O}(hy+\lambda_1^{-y})={\rm O}(hy)=$$
$$
{\rm O}\left(h\frac1{\ln\lambda_1} \ln\frac1{|h|}\right)={\rm O}\left(h\ln\frac1{|h|}\right)$$ for $0\leq h \leq e^{-\frac12}$, hence (\ref{regFm}) holds.

Now, for $0 < \kappa < 1$ and $0 < h \leq 1$, we have $h \, \ln {1 \over h} \leq {1 \over 1 - \kappa} \, h^\kappa$.\footnote{This inequality is equivalent to $\ln \frac1{h^{1-\kappa}}\leq \frac1{h^{1-\kappa}}$.} Therefore $F$ is H\"olderian with exponent $\kappa$.

3) Let $G(x) = \sum_{k \geq 1} {1\over u_k} \, \sin (2\pi  v_k x)$, where $(v_k)$ is a sequence of integers
and $(u_k)$ is a sequence of positive numbers such that $\sum 1/ u_k < \infty$. If $(t_k)$ is an increasing sequence of integers,
the ergodic sums of $G$ at time $t_n$ reads:
\begin{eqnarray}
G^{(t_n)}(x) =  \sum_{k \geq 1} \, {1 \over u_k} \, {\sin(\pi  v_k t_n \alpha) \over \sin(\pi  v_k \alpha)}
\sin(\pi  v_k (2x + (t_n -1) \alpha)) = A_n + B_n + C_n, \label{sumerg1}
\end{eqnarray}
where in (\ref{sumerg1}) $A_n, B_n, C_n $ are respectively the partial sums $\sum_{k < n}$, $\sum_{k = n}$, $\sum_{k > n}$.

Now, take $G = F_s$ given by (\ref{defF1}) and consider $t_n = q_{m_n}$. The decomposition (\ref{sumerg1}) yields (for some constant $c'>0$):
\begin{equation}
|A_n| \leq \sum_{k \leq n-1} \, {1 \over q_{m_k}} \, {|\sin(\pi  q_{m_k} q_{m_n} \alpha)| \over|\sin(\pi q_{m_k} \alpha)|}
\leq c' {1\over  q_{m_n+1}} \, \sum_{k \leq n-1} \,  q_{m_k};\end{equation}
using~(\ref{minbpq1}),
\begin{equation}
|B_n|  = {1 \over q_{m_n}} \, {|\sin(\pi  q_{m_n} q_{m_n} \alpha)| \over|\sin(\pi q_{m_n} \alpha)|}  \asymp \, 1;
\end{equation}
\begin{equation}
|C_n| \leq \sum_{k \geq n+1} \, {1 \over q_{m_k}} \, {|\sin(\pi  q_{m_k} q_{m_n} \alpha)| \over|\sin(\pi q_{m_k} \alpha)|}
\leq c' \, q_{m_n} \sum_{k \geq n+1} \, {1\over q_{m_k}}.
\end{equation}
By (\ref{maj1qn}) and~(\ref{maj2qn}), we have
$$|A_n| \leq 2 {q_{m_{n-1}+1} \over  q_{m_n+1}}, \ |C_n| \leq C {q_{m_{n}} \over  q_{m_{n+1}}}.$$

It follows that we can select the sequence $s = (m_k)$ such that the terms $A_n, C_n$ cannot cancel the contribution of $B_n$. That is, $B_n$ is bounded away from zero (and is clearly bounded) and the behavior of $F_s^{(q_{m_n})}$ is similar to the behavior of $B_n\sin(2\pi q_{m_n}x)$. This, by Lemma~\ref{supp} yields an uncountable set of essential values of $F_s$.

Remark that we have boundedness of the ergodic sums of $F_1$ at time $q_n$: $\|F_1^{(q_n)}\|_2 \leq 2\pi$.
It follows that, for $\delta > 0$ small enough, the above property of the existence of an uncountable set of essential values for $F_s$ is still satisfied
for $F_s + \delta F_1$ (in other words, we obtain a stability of ergodicity of $F_s$ by some perturbations).
\eop

\begin{center}
\includegraphics[scale=.12]{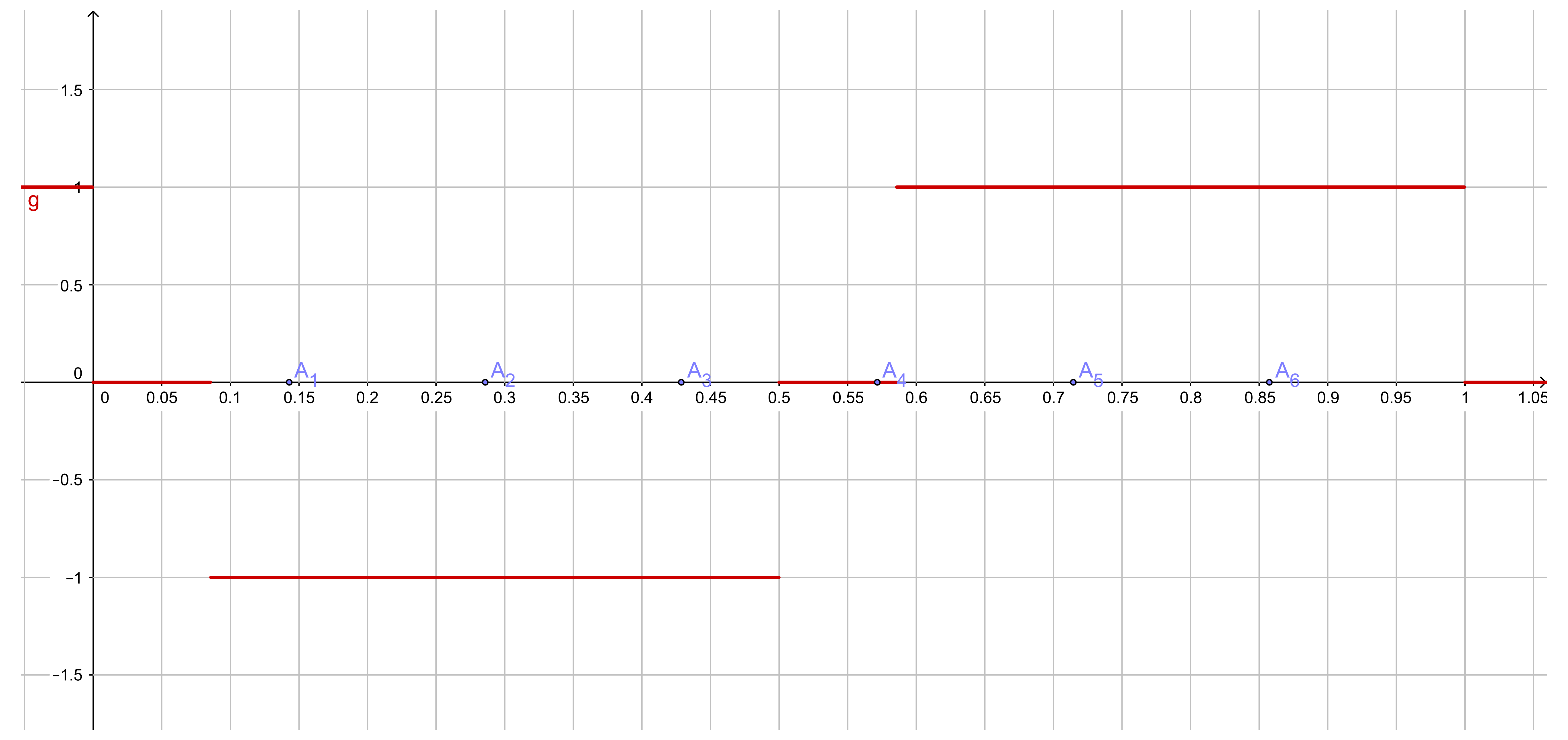} \\
{\it Fig.1 \ rotation $\alpha=\pi -3$, $\beta= 2- \sqrt 2$, $\Phi = \Phi_\beta$,  $A_k= {k\over7}$, graph of $\Phi_\beta$}

\includegraphics[scale=.12]{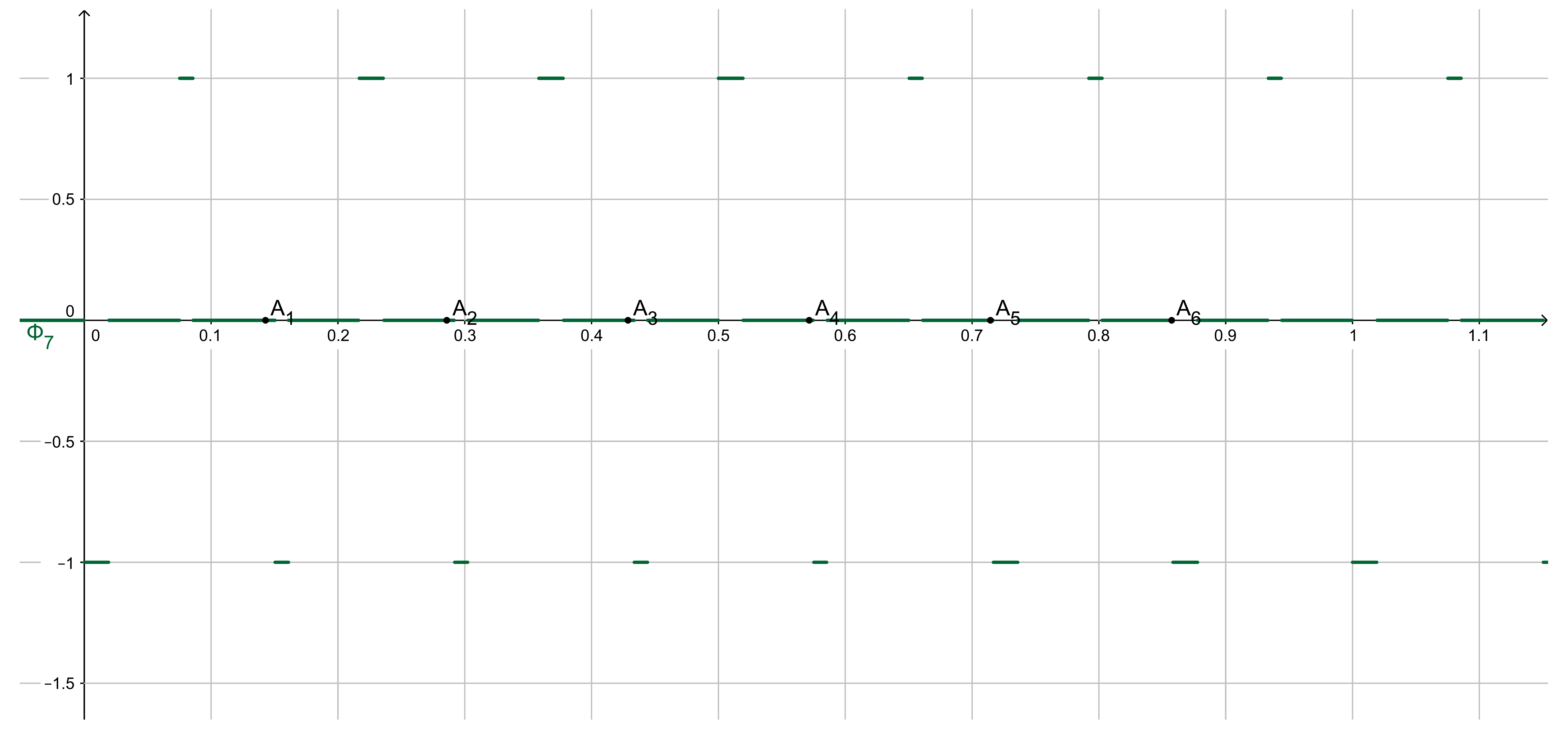} \\
{\it Fig.2  \ $\Phi_7 = \Phi_\beta^{(7)}$,  $A_k= {k\over7}$, graph of $ \Phi_\beta^{(7)}$}

\includegraphics[scale=0.12]{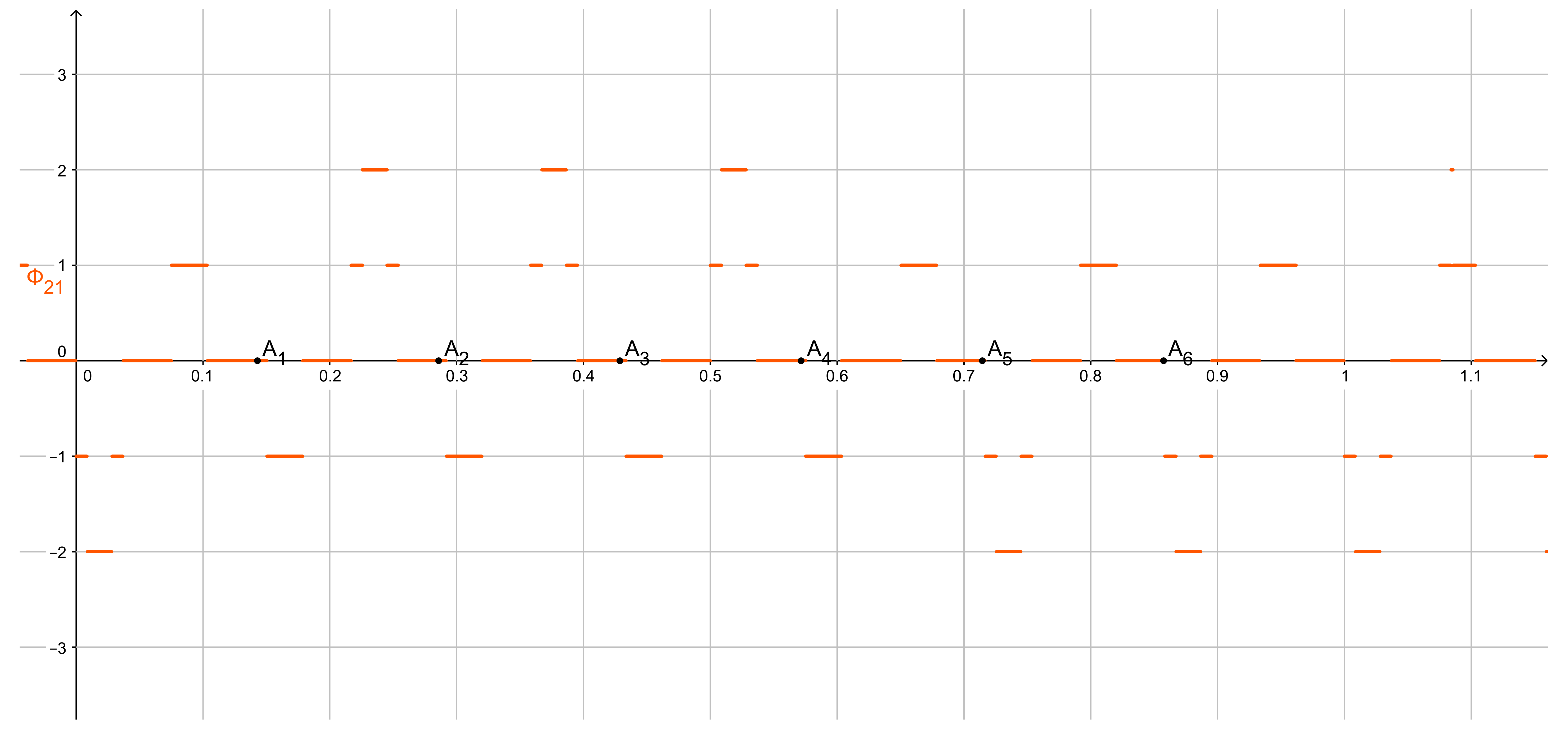} \\
{\it Fig.3  \ $\Phi_{21} = \Phi_\beta^{(21)} = \Phi_\beta^{(7)} + \Phi_\beta^{(7)}(. + 7\alpha) +\Phi_\beta^{(7)}(. + 14 \alpha)$, $A_k= {k\over7}$, graph of $\Phi_\beta^{(21)}$}
\end{center}

\vskip 6mm
\end{document}